\documentclass[preprint,3p,12pt]{elsarticle}
\usepackage[utf8]{inputenc}
\usepackage{textcomp}

\usepackage{graphicx}
\usepackage{amsmath,amssymb}
\usepackage[version=4]{mhchem}
\usepackage{siunitx}
\usepackage{longtable,tabularx}
\usepackage[symbol]{footmisc}

\usepackage{bm}
\usepackage{amsfonts}
\usepackage{amscd}
\usepackage{hyperref}
\usepackage{pstricks}
\usepackage{subcaption}
\usepackage{framed,fancybox}
\usepackage{bbm}
\usepackage{multirow}
\usepackage{threeparttable}
\usepackage{rotating}
\usepackage{stfloats}
\usepackage{color,soul}
\usepackage{float}
\usepackage{afterpage}
\usepackage{listings}
\usepackage{xcolor}

\usepackage{amsthm}

\definecolor{codegreen}{rgb}{0,0.6,0}
\definecolor{codegray}{rgb}{0.5,0.5,0.5}

\lstdefinestyle{mystyle}{
  commentstyle=\color{codegreen},
  numberstyle=\tiny\color{codegray},
  keywordstyle=\color{blue},
  basicstyle=\ttfamily\small,
  breakatwhitespace=false,         
  breaklines=true,                 
  keepspaces=true,                 
  numbers=left,                    
  numbersep=5pt,                  
  showspaces=false,                
  showstringspaces=false,
  showtabs=false,                  
  tabsize=2
}

\allowdisplaybreaks

\journal{Journal of the Franklin Institute}

\begin{document}

\begin{frontmatter}

\title{{\bf Optimal Control of Parabolic Differential} \\{\bf  Equations Using Radau Collocation}\footnote[2]{Portions of this work were presented at the 2024 and 2025 American Control Conferences (ACC) in Toronto, Ontario, Canada, and Denver, Colorado (https://doi.org/10.23919/ACC60939.2024.10644743, https://doi.org/10.23919/ACC63710.2025.11107475). Portions of the work were also presented at the 2025 AAS/AIAA Spaceflight Mechanics Meeting in Lihue, Hawaii.}}

\author[ufaffil]{Alexander M.~Davies\fnref{label1}}
\ead{alexanderdavies@ufl.edu}
\fntext[label1]{Ph.D. Candidate, Department of Mechanical and Aerospace Engineering.}
\affiliation[ufaffil]{organization={University of Florida},
            city={Gainesville},
            postcode={32611},
            state={FL},
            country={USA}}

\affiliation[afrlaffil]{organization={Air Force Research Laboratory},
            city={Eglin AFB},
            postcode={32542},
            state={FL},
            country={USA}}
            
\author[ufaffil]{Sara Pollock\fnref{label2}}
\ead{s.pollock@ufl.edu}
\fntext[label2]{Associated Professor, Department of Mathematics.}

\author[afrlaffil]{Miriam E.~Dennis\fnref{label3}}
\ead{miriam.dennis.1@us.af.mil}
\fntext[label3]{Research Engineer, Munitions Directorate.}

\author[ufaffil]{Anil V. Rao\corref{cor1}\fnref{label4}}
\ead{anilvrao@ufl.edu}
\cortext[cor1]{Corresponding author.}
\fntext[label4]{Professor, Department of Mechanical and Aerospace Engineering. AIAA Associate Fellow. AAS Fellow.}

\begin{abstract}
A method is presented for the numerical solution of optimal boundary control problems governed by parabolic partial differential equations. The continuous space-time optimal control problem is transcribed into a sparse nonlinear programming problem through state and control parameterization. In particular, a multi-interval flipped Legendre-Gauss-Radau collocation method is implemented for temporal discretization alongside a Galerkin finite element spatial discretization. The finite element discretization allows for a reduction in problem size and avoids the redefinition of constraints required under a previous method. Further, a generalization of a Kirchoff transformation is performed to handle variational form nonlinearities in the context of numerical optimization. Due to the correspondence between the collocation points and the applied boundary conditions, the multi-interval flipped Legendre-Gauss-Radau collocation method is demonstrated to be preferable over the standard Legendre-Gauss-Radau collocation method for optimal control problems governed by parabolic partial differential equations. The details of the resulting transcription of the optimal control problem into a nonlinear programming problem are provided. Numerical examples demonstrate that the use of a multi-interval flipped Legendre-Gauss-Radau temporal discretization can lead to a reduction in the required number of collocation points to compute accurate values of the optimal objective in comparison to other methods. Lastly, a self-convergence analysis on each test problem illustrates that the error decays exponentially as a function of the mesh size in both the temporal and spatial dimensions.
\end{abstract}

\begin{keyword}
    parabolic partial differential equations \sep Gaussian quadrature collocation \sep finite element methods \sep optimal control \sep nonlinear optimization
\end{keyword}

\end{frontmatter}

\section{Introduction}
The study of optimal control has largely focused on systems governed by ordinary differential equations (ODEs) \cite{Rao2009}; however, in applications where spatial effects cannot be ignored, higher-fidelity models are better suited to describe the dynamical behavior. To be more specific, the optimal control of systems governed by partial differential equations (PDEs) \cite{Lions1971,Troltzsch2024,ManzoniQuarteroni2021,BeckerKapp2000} has garnered more interest in recent decades as computational resources continue to grow. Some fields in which the optimal control of PDEs has been applied include, but are not limited to:  fluid mechanics \cite{Casas1998,Ghattas1997,AllahverdiPozo2016,Sritharan1998,holl2020,RabaultKutcha2019,MortonJameson2018,HwangLee2022}, thermodynamics \cite{DebusscheFuhrman2007,CasasTroltzsch2020,CleverLang2012}, finance \cite{KafashDelavarkhalafi2016}, and epidemiology \cite{WangZhang2023}. 
Generally, PDE-constrained optimal control problems do not possess analytical solutions, which implies that efficient numerical methods are required for their solution.

PDE-constrained optimal control problems present a few nuances in comparison to ODE-constrained optimal control. For instance, the added dimensionality of a PDE system promotes an increased emphasis on the accuracy and efficiency of the applied discretization. As the number of dimensions grows, the size of the problem follows suit, which implies that poorly-selected discretizations can cause a rapid increase in the problem size in comparison to ODE-constrained problems. Further, the addition of a spatial variable implies the introduction of \textit{boundary conditions} that may or may not be a function of the control variable. Due to the presence of the boundary conditions, it is imperative that the selected discretization is well-posed to constrain the apparent boundary conditions, which is not a consideration for ODE systems. Lastly, the infinite-dimensional character of PDE solution spaces amplifies the importance of the quality of the discrete approximation. Identifying an appropriate finite-dimensional subspace in which to approximate the solution for both the temporal and spatial dimension becomes a central consideration for numerical methods for PDE-constrained optimal control.


Numerical methods for optimal control often fall under two categories: indirect and direct methods. Indirect methods transform the original optimal control problem into a Hamiltonian boundary value problem via the derivation of the first-order optimality conditions from the calculus of variations. The Hamiltonian boundary value problem is then solved by a shooting, multiple-shooting, or collocation approach \cite{Kirk2004}. Direct methods transcribe the original optimal control problem into a large, sparse nonlinear programming problem (NLP) through parameterization of the control and/or state.  The NLP can then be solved with well-developed NLP solvers such as \textit{SNOPT} \cite{GillMurray2002}, \textit{IPOPT} \cite{BieglerZavala2008}, or \textit{KNITRO} \cite{ByrdNocedal2006}. 

A benefit of the indirect approach is that one can verify that the solution obtained is optimal (for the selected discretization) given that the optimality conditions are satisfied. A drawback, however, is that the first-order optimality conditions are \textit{problem-specific} and can be complex to derive when the underlying PDE is highly nonlinear. In addition, good initial guesses are required for the state, adjoint, and control parameterizations, which can be difficult if no analytical solution is available or if the dynamics are particularly stiff. In addition, for problems with state-inequality constraints, the optimality conditions are further complicated if the state lies on a constrained arc \cite{Betts2020, ByczkowskiRao2024}, as the optimality conditions change over active regions. This would motivate a multi-phase approach, which can be difficult to implement  without a priori knowledge of the time at which the constrained arc becomes active. 

Direct methods, at the expense of \textit{verified} optimality (though there are instances in which the direct and indirect approaches have been shown to coincide \cite{BeckerVexler2007}), mitigate some of the drawbacks of the indirect approach. Direct methods remove the adjoint variable from the problem formulation entirely, which reduces the complexity of discretizing the optimal control problem. Further, direct methods present a larger domain of convergence than do indirect methods, which implies that a good initial guess is not necessarily required when solving the problem with a direct method. Also, due to the versatility of the direct transcription approach, inequality constraints can be easily enforced without a priori knowledge of constrained arcs. Direct methods are a highly general approach that can be applied to a wide array of optimal control problems constrained by parabolic systems, and as a result, we are motivated to employ a direct method in this paper. 

Historically, discretization methods for the numerical solution of optimal control problems constrained by PDEs have been classified as either $h-$methods or $p-$methods. In an $h-$method, typically, a low-order polynomial is used to approximate the state, and the domain is divided into a set of intervals. Convergence is achieved by increasing the number of intervals. Some examples of $h-$methods that have been applied to PDE-constrained optimal control problems are Euler methods and Runge-Kutta methods \cite{Betts2020, BuskensGriesse2006,NeitzelTroltzsch2009, BurgerMachbub1991}. In a $p-$method, a single interval is used, and the degree of the state approximation is increased to achieve convergence. Examples of $p-$methods include spectral and pseudospectral methods \cite{AliShamsi2019, KhaksarShamsi2017, ChenHuang2017}. The combination of both $h-$ and $p-$ methods is known as an $hp-$method, where both the number of intervals and the degree of the state approximation in each interval can be increased to achieve convergence.

One efficient class of numerical methods that has grown increasingly popular for the optimal control of ODE systems is {\em direct orthogonal collocation} \cite{BensonHuntington2006,VlassenbroeckVanDooren1988,Williams2004,GargPatterson2011,GargPatterson2010}. In an orthogonal collocation method, the support points of the state and control approximations are determined by the roots of orthogonal polynomials such as the Chebyshev, Legendre, or Jacobi polynomials. Sometimes, these methods are referred to as {\em Gaussian quadrature collocation} methods, as the support points are often associated with Gaussian quadrature rules. Gaussian quadrature orthogonal collocation methods can be cast within an $hp-$framework, which is a primary focus of this work. This formulation provides advantages over existing approaches by enabling multiple avenues to achieve convergence. When the solutions to optimal control problems are smooth, it has been demonstrated in Ref.~\cite{HagerHou2019} that Gaussian quadrature collocation methods converge to highly-accurate solutions at an exponential rate. 

Gaussian quadrature orthogonal collocation methods have been scarcely applied to 
the optimal control of PDEs. 
The few methods that have been developed have been formulated as $p-$methods \cite{KhaksarShamsi2017,LiZhou2018,AliShamsi2019}, where the state is approximated by a single (or global) basis of polynomials. In Ref.~\cite{DarbyHager2011a}, it was demonstrated that a {\em multi-interval} (or $hp-$) orthogonal collocation approach can increase the utility of Gaussian quadrature orthogonal collocation while remaining close to its desirable exponential convergence properties. The work of Ref.~\cite{DaviesDennis2024} develops a multi-interval orthogonal collocation approach (with its application in both space and time) for the solution of optimal control problems governed by time-dependent PDEs to handle parabolic problems and general sets of boundary conditions. However, the use of a modified Legendre-Gauss-Radau (mLGR) \cite{EideHager2021} set of constraints was required to satisfy required state relationships and Neumann boundary conditions. 
This, in turn, led to an overconstrained optimization problem necessitating the redefinition of the mLGR constraints as inequality constraints between a user-specified tolerance, and as a result, an avenue was introduced for increased error. The method developed in this paper serves as an improvement to the method developed in Ref.~\cite{DaviesDennis2024}, where, here, a finite element discretization is employed in the spatial dimension to avoid the disadvantages of the spatial, multi-interval, orthogonal collocation approach employed in Ref.~\cite{DaviesDennis2024}.

Finite element discretizations have been widely applied to optimal control problems governed by PDEs \cite{BurgerPogu1991,KupferSachs1992,Heinkenschloss1996,BurgerMachbub1991,MeidnerVexler2008a,BeckerMeidner2007,BeckerVexler2007}. The finite element method discretizes the variational form of the partial differential equation which requires that integral averages be satisfied as opposed to satisfying pointwise equalities.  Namely, the variational form is constructed through the multiplication of the PDE by a test function and integration over the spatial domain. Typically, the solution is assumed to belong to a Hilbert space, and a finite set of basis functions are utilized to approximate the solution structure. Nonlinearities in the variational form of a PDE can lead to complications in forming a satisfactory discretization of the equation. In the absence of an optimization routine, several well-known techniques exist to address nonlinear problems. In general, in the framework of direct collocation methods with finite elements for the solution of optimal control problems governed by PDEs, few methods exist to handle nonlinear terms in the variational formulation. 

Though orthogonal collocation methods and finite element approaches have been used extensively in their respective application areas, no work has been done to combine the approaches in a unified, generalized framework for the direct solution of optimal control problems governed by parabolic PDEs. Generally, the focus of orthogonal collocation has been primarily applied to optimal control problems governed by ODEs. For classes of problems governed by ODEs, it has been shown that methods of Radau-type and flipped Radau-type may be used essentially interchangeably to achieve solutions of similar accuracy \cite{GargPatterson2010}. In this work, it is shown that this is \textit{not} the case for problems governed by parabolic PDEs, and that methods of flipped-Radau type are superior to standard Radau-type methods for this particular application due to the well-posedness of the discretization. Extensive research has been dedicated to the application of the finite element method to optimal control problems governed by PDEs, especially in the context of indirect methods. However, little work has been done to address nonlinearities in the variational formulation within the context of direct collocation methods. This work seeks to mitigate the computational burden that stems from nonlinearities in the variational framework through a ``Kirchoff-like'' integral linearization. The approximation aids to reduce the computational complexity of evaluating and computing NLP derivatives of nonlinear nonlinear-type terms within the optimization context.  

Ultimately, the novelty and contribution of this work is as follows. A novel method is presented for the numerical solution of optimal control problems governed by parabolic partial differential equations. A multi-interval flipped Legendre-Gauss-Radau (fLGR) temporal discretization is implemented alongside a Galerkin finite element spatial discretization for the direct solution of optimal control problems governed by PDEs. The finite element discretization serves as an improvement to a previous method \cite{DaviesDennis2024} through maintaining a proper NLP search space and the removal of the requirement for the redefinition of constraints to handle boundary conditions. In addition, the finite element discretization and the variational approach allow for a generalization to higher dimensions with varying domain shapes. Furthermore, the location of the collocation points under the fLGR temporal discretization allows for the satisfaction of dynamics/boundary relationships at the terminal time, which is not possible for explicit methods in an optimization framework with boundary control. Another contribution is the application of a ``Kirchoff-like'' transformation to nonlinear terms in time-dependent PDEs within the optimization framework. 
The transformation allows for the direct solution of a larger, linear system with additional constraints imposed with the optimization framework. Further, numerical comparisons are performed to highlight the effectiveness of the novel method in the context of existing solutions in the literature. The use of the novel approach is shown to be capable of computing solutions with similar values of the optimal objective to other methods in fewer numbers of points.  

The remainder of this paper is structured as follows. The notations and conventions used in this work are provided in Section \ref{sec:notcov}. A general optimal control problem governed by a nonlinear 1-D parabolic PDE is described and transformed into a multi-interval formulation in Section \ref{sec:OCP}. The variational form of the PDE is outlined and the finite element discretization is implemented in Section \ref{sec:FEM}. The fLGR temporal discretization is illustrated and applied in Section \ref{sec:fLGR}. The transcription of the discretized problem into an NLP (as well as NLP derivatives, sparsity patterns, etc.) is outlined in Section \ref{sec:NLP}. The method is then utilized on a viscous Burgers' tracking problem \cite{BuskensGriesse2006,Betts2020} and a nonlinear heat equation problem \cite{Heinkenschloss1996,Betts2020} and compared with existing solutions in Section \ref{sec:NumericalExamples}. Finally, Section \ref{sec:conclusions} provides conclusions on this research. 

\section{Notation and Conventions}\label{sec:notcov}
In this paper, some notational shortcuts are used to aid the reader and maintain conciseness. Relevant notation and conventions are provided here for reference. All vectors are denoted by a bold marking. Matrices are not bolded. Given the vector $\mathbf{f} \in \mathbb{R}^{1 \times n}$ and vector $\mathbf{g} \in \mathbb{R}^{1 \times n}$, an example of a relevant operation would be 
$$
    F = \mathbf{f}^T\mathbf{g} \quad \in \mathbb{R}^{n \times n},
$$
where $\mathbf{f}^T \in \mathbb{R}^{n \times 1}$ denotes the transpose of vector $\mathbf{f}$. Likewise to matrices, scalars are not bolded, a similar operation provides
$$    
f = \mathbf{f}\mathbf{g}^T \quad \in \mathbb{R}. 
$$
To refer to a specific element of a vector, subscript notation is used. That is, if the vector $\mathbf{f} \in \mathbb{R}^{1 \times n}$ is defined by 
$$
    \mathbf{f} = \begin{bmatrix}
        f_1 & f_2 & \ldots & f_n
    \end{bmatrix},
$$
the $i^{\mathrm{th}}$ element of vector $\mathbf{f}$ is denoted by $f_i$. For matrices in this work, a single subscript is used to indicate a specific column of the matrix. Given the matrix $F \in \mathbb{R}^{n \times n}$ written as
$$
    F  = \begin{bmatrix}
        \mathbf{F}_1 & \mathbf{F}_2 & \ldots & \mathbf{F}_n
    \end{bmatrix},
$$
the $i^{\mathrm{th}}$ column of matrix $F$ is denoted by $\mathbf{F}_i \in \mathbb{R}^{n \times 1}$. When referring to a specific column of a matrix, the matrix variable is bolded with a subscript to refer to a specific column of that matrix. To refer to a specific subset of indices, a complete list of notation is provided:
\begin{itemize}
    \item Given the matrix $F \in \mathbb{R}^{n \times n}$, one may extract a specific subset of columns by 
    $$
    F_{i:j} = \begin{bmatrix}
        \mathbf{F}_i & \ldots & \mathbf{F}_j
    \end{bmatrix},
    $$
    where $j \geq i$.
    \item Given the matrix $F \in \mathbb{R}^{n \times n}$, one may extract a specific subset of rows by 
    $$
    F_{i:j,:} \triangleq \mathrm{all\:of\:the\:columns\:in\:the\:}i^{\mathrm{th}}\:\mathrm{through}\:j^{\mathrm{th}}\:\mathrm{row\:of}\:F,
    $$
    where $j \geq i$.
    \item Given the matrix $F \in \mathbb{R}^{n \times n}$, to extract one column of $F$, the result is a vector, thus, we have
    $$
    \mathbf{F}_i \triangleq i^{\mathrm{th}}\:\mathrm{column}\:\mathrm{of}\:F.
    $$
    \item Given the matrix $F \in \mathbb{R}^{n \times n}$, to refer to a specific element of $F$, we use double subscripts by
    $$
    F_{ij} \triangleq \mathrm{element\:of\:matrix\:}F\:\mathrm{in\:the\:}i^{\mathrm{th}}\:\mathrm{row\:and}\:j^{\mathrm{th}}\:\mathrm{column}
    $$
    \item Given the vector $\mathbf{f} \in \mathbb{R}^{1 \times n}$, one may extract a specific subset of columns by 
    $$
    \mathbf{f}_{i:j}   = \begin{bmatrix}
        f_i & \ldots & f_j
    \end{bmatrix},
    $$
     where $j \geq i$.
     \item Given the vector $\mathbf{f} \in \mathbb{R}^{1 \times n}$, to extract one element of $\mathbf{f}$, the result is a scalar, thus, we have
     $$
     f_i \triangleq i^{\mathrm{th}}\:\mathrm{element}\:\mathrm{of}\:\mathbf{f}.
     $$
\end{itemize}
Further, some important vector and matrix operations are defined as follows. Given the vector $\mathbf{f}(x,t) \in \mathbb{R}^{1 \times n}$, the partial derivative with respect to $x$ is defined as 
$$
    \dfrac{\partial \mathbf{f}(x,t)}{\partial x} = \begin{bmatrix}
        \dfrac{\partial f_1(x,t)}{\partial x} & \dfrac{\partial f_2(x,t)}{\partial x} & \ldots & \dfrac{\partial f_n(x,t)}{\partial x}
    \end{bmatrix}.
$$
The gradient (or Jacobian) of $\mathbf{f}^T$ with respect to a vector $\mathbf{z} \in \mathbb{R}^{1 \times n}$ is defined by 
$$
\nabla_{\mathbf{z}}\mathbf{f}^T = \begin{bmatrix}
    \frac{\partial f_1}{\partial z_1} & \frac{\partial f_1}{\partial z_2} & \ldots & \frac{\partial f_1}{\partial z_n} \\
    \frac{\partial f_2}{\partial z_1} & \frac{\partial f_2}{\partial z_2} & \ldots & \frac{\partial f_2}{\partial z_n} \\
    \vdots & \vdots & \ddots & \vdots\\
    \frac{\partial f_n}{\partial z_1} & \frac{\partial f_n}{\partial z_2} & \ldots & \frac{\partial f_n}{\partial z_n}
\end{bmatrix}.
$$
To strip the elements of $\nabla_{\mathbf{z}}\mathbf{f}^T$ into a single column, the notation 
$$
\nabla_{\mathbf{z}}\mathbf{f}^T (:) = \begin{bmatrix}
    \left(\frac{\partial \mathbf{f}^T}{\partial z_1}\right)^T & \left(\frac{\partial \mathbf{f}^T}{\partial z_2}\right)^T & \ldots &\left(\frac{\partial \mathbf{f}^T}{\partial z_n}\right)^T
\end{bmatrix}^T,
$$
is used. Additional detail and explanation of the notation will be provided as required.
\section{Optimal Control of a Parabolic Partial Differential Equation}\label{sec:OCP}
Without loss of generality, an optimal control problem is presented in Lagrange form as follows. The goal is to minimize
\begin{equation}
    \mathcal{J} = \int_{t_0}^{t_f}\int_{\Omega} \mathcal{L}(x,t,y(x,t))\:\mathrm{d}x\:\mathrm{d}t+\int_{t_0}^{t_f} \mathcal{P}(t,u_1(t),u_2(t))\:\mathrm{d}t, \label{eq:objective}
\end{equation}
subject to the one-dimensional nonlinear PDE  
\begin{equation}
    \frac{\partial y}{\partial t} + \kappa(y) \frac{\partial y}{\partial x} = a \frac{\partial ^2 y}{\partial x^2}, 
\end{equation}
with boundary controls of Neumann-kind and fixed-state initial conditions
\begin{align}
    \frac{\partial y(0,t)}{\partial x} &= u_1(t), \label{eq:bc1} \\
    \frac{\partial y(1,t)}{\partial x} &= u_2(t), \label{eq:bc2} \\
     y(x,t_0) &= q(x,t_0), \label{eq:initialcondition}
\end{align}
where $\kappa(y)$ is an integrable function of the state variable, 
$q(x,t_0)$ is a specified set of initial conditions, and $a$ is a constant parameter. The problem is defined on $t_0 \leq t \leq t_f$ and $\Omega \in [0,1]$. The temporal mesh is mapped onto the fixed interval $\tau \in [-1,+1]$ by 
\begin{equation}
    \tau = \frac{2(t-t_0)}{t_f-t_0}-1,
\end{equation}
which leads to the transformed PDE 
\begin{equation}
    \frac{\partial y}{\partial \tau} + \frac{\partial t}{\partial \tau} \kappa(y) \frac{\partial y}{\partial x} = a \frac{\partial t}{\partial \tau}\frac{\partial ^2 y}{\partial x^2}, \label{PDE1}
\end{equation}
where ($\partial t/\partial \tau$) is given by 
\begin{equation}
    \frac{\partial t}{\partial \tau} = \frac{t_f-t_0}{2}.
\end{equation}
The mesh is subdivided into $J$ time intervals, $T_j \in [\tau_{j-1},\tau_j]$, $j = \{1,\ldots,J\}$, where ($\tau_0,\ldots,\tau_J$) are the mesh points. The mesh intervals have the properties that $\cup_{j=1}^J T_j = [-1,1]$ and $\cap_{j=1}^J T_j = 0$, while the mesh points are such that $-1 = \tau_0 < \tau_1 <\tau_2< \ldots < \tau_J = 1.$ Each interval, $T_j$, is mapped to $r \in [-1,1]$ by the transformations
\begin{equation}
    r = 2\frac{\tau-\tau_{j-1}}{\tau_j-\tau_{j-1}}-1,\quad \tau \in [\tau_{j-1},\tau_j], \quad j = \{1,\ldots,J\}.
\end{equation}
The affine transformations yield the transformed PDE 
\begin{equation}
     \frac{\partial y^{(j)}}{\partial r} +  \frac{\partial t}{\partial \tau}\left(\frac{\partial \tau}{\partial r} \right)^{(j)} \kappa(y^{(j)}) \frac{\partial y^{(j)}}{\partial x} = a \frac{\partial t}{\partial \tau}\left(\frac{\partial \tau}{\partial r} \right)^{(j)} \frac{\partial^2 y^{(j)}}{\partial x^2}, 
\end{equation}
with the interval mappings, $(\partial \tau/\partial r)^{(j)}$, defined by 
\begin{equation}
    \left(\frac{\partial \tau}{\partial r} \right)^{(j)} = \frac{(\tau_j-\tau_{j-1})}{2}, \quad j = \{1,\ldots J\}.
\end{equation}
For the sake of brevity in the remaining discussion, the mappings are grouped by 
\begin{equation}
    \psi^{(j)} = \frac{\partial t}{\partial \tau}\left(\frac{\partial \tau}{\partial r} \right)^{(j)}\quad j = \{1,\ldots J\}.
\end{equation}
Thus, the transformed PDE can ultimately be rewritten as
\begin{equation}
     \frac{\partial y^{(j)}}{\partial r} +  \psi^{(j)} \kappa(y^{(j)}) \frac{\partial y^{(j)}}{\partial x} = a \psi^{(j)} \frac{\partial^2 y^{(j)}}{\partial x^2}, \quad j = \{1,\ldots J\}. \label{eq:strongpde}
\end{equation}
The PDE is now enforced to hold locally over each temporal interval, and additional constraints are required to enforce 
equality of the state at interval boundaries. 
This approach is addressed in Section \ref{sec:fLGR}.  

\section{The Variational Form and the Finite Element Discretization}\label{sec:FEM}
 A procedure is now described for defining a finite element discretization in space.  The method described in this section is preferred because of its flexibility and ability to generalize to higher order methods. However, a finite difference method can also be used in its place, the details of which appear in the Appendix. 

A variational form of the transformed PDE can be constructed through the multiplication of the equation by a test function from a suitable space and integration over the spatial domain. The variational form of the PDE is provided by 
\begin{equation}
\int_{\Omega}  \frac{\partial y^{(j)}}{\partial r} v^{(j)}\:\mathrm{d}x + \psi^{(j)} \int_{\Omega} \kappa(y^{(j)}) \frac{\partial y^{(j)}}{\partial x} v^{(j)}\:\mathrm{d}x = a \psi^{(j)} \int_{\Omega} \frac{\partial^2 y^{(j)}}{\partial x^2} v^{(j)}\:\mathrm{d}x, \label{eq:variational1}
\end{equation}
where $v^{(j)}$ is the test function in 
interval $T_j$, 
$j=\{1,\ldots,J\}$, defined in $V = \{v^{(j)} \in H_1(\Omega) \::\: v^{(j)}\in L_2(\Omega),\: \frac{\partial v^{(j)}}{\partial x} \in L_2(\Omega)\}$ with $L_2(\Omega)$ given by the space of square-integrable functions on $\Omega$ such that $ \forall f \in L_2(\Omega)$:
\begin{equation}
    ||f||_{L_2(\Omega)} \triangleq \left( \int_{\Omega} |f|^2\:\mathrm{d}x \right)^{\frac{1}{2}} < \infty.
\end{equation}
Integrating Eq.~\eqref{eq:variational1} by parts provides
\begin{multline}
\int_{\Omega}  \frac{\partial y^{(j)}}{\partial r} v^{(j)}\:\mathrm{d}x + \psi^{(j)} \int_{\Omega} \kappa(y^{(j)}) \frac{\partial y^{(j)}}{\partial x} v^{(j)}\:\mathrm{d}x  \\= a \psi^{(j)} \left(\left.\frac{\partial y^{(j)}}{\partial x}v^{(j)}\right|_{x=1} - \left.\frac{\partial y^{(j)}}{\partial x}v^{(j)}\right|_{x=0} - \int_{\Omega} \frac{\partial y^{(j)}}{\partial x} \frac{\partial v^{(j)}}{\partial x}\:\mathrm{d}x \right),
\end{multline}
and substituting the relationships for the boundary conditions gives 
\begin{multline}
\int_{\Omega}  \frac{\partial y^{(j)}}{\partial r} v^{(j)}\:\mathrm{d}x + \psi^{(j)} \int_{\Omega} \kappa(y^{(j)}) \frac{\partial y^{(j)}}{\partial x} v^{(j)}\:\mathrm{d}x \\= a \psi^{(j)} \left(u_2^{(j)}v^{(j)}(1) - u_1^{(j)}v^{(j)}(0)
- \int_{\Omega} \frac{\partial y^{(j)}}{\partial x} \frac{\partial v^{(j)}}{\partial x}\:\mathrm{d}x \right). \label{eq:weakform1}
\end{multline}

Thus, the new problem is to find a solution $y^{(j)}$, $j = \{1,\ldots,J\}$, that satisfies Eq.~\eqref{eq:weakform1} for all test functions $v^{(j)}$, $j = \{1,\ldots,J\}$, in the test space $V$. In this work, the Galerkin approximation is assumed, where the space in which the solution is to be sought is identical to the test space $V$. The first step under the Galerkin finite element method is to approximate the infinite-dimensional form of the PDE in Eq.~\eqref{eq:weakform1} through finite-dimensional functions, $y_h$ and $v_h$, that are constructed by linear combinations of basis functions that define a finite-dimensional subspace of $H_1(\Omega)$. These approximations can be written by 
\begin{align}
    y^{(j)} &\approx y^{(j)}_h  = \sum_{k = 1}^{N_x} \phi_k Y^{(j)}_k, \label{findim1} \\
    v^{(j)} &\approx v^{(j)}_h   = \sum_{i = 1}^{N_x} \phi_i V^{(j)}_i, \label{findim2}
\end{align}
where $N_x$ represents the number of nodes or basis functions utilized to approximate the functions $y^{(j)}$ and $v^{(j)}$. Eventually, the finite-dimensional approximation of the PDE is transcribed to a discrete form to be handled by a computer. 
\subsection{Nonlinearities in the Optimization Framework}
As the solution procedure is implicit, the nonlinear term may be computed through direct integration of the nonlinear form in Eq.~(\ref{eq:weakform1}) upon every NLP iteration. However, this makes evaluation of the NLP constraints function rather computationally expensive. Further, as many NLP solvers require derivatives of the NLP functions, analytical derivation of the NLP Jacobian of the constraints becomes a far more arduous task. One can seek to mitigate this difficulty through automatic differentiation software; however, computational speed is often limited by coding restrictions and file sizes that scale with the mesh. As the mesh grows larger, evaluation of the nonlinear term grows more expensive, derivatives of the constraints are likewise more expensive, and generated file sizes can grow out of proportion with the size of the mesh. 

Due to the drawbacks of handling the term directly, we seek an alternate procedure to resolve, or in this case, approximate the nonlinearity. One possible technique is to rewrite Eq.~\eqref{eq:weakform1} via a ``Kirchoff-like'' integral transformation. Kirchoff transformations are often employed 
in applications to the heat equation when the thermal conductivity depends upon the temperature. The standard Kirchoff transformation may be expressed as follows. Consider the 1-D steady-state nonlinear heat conduction equation:
\begin{equation}
    \dfrac{\mathrm{d}}{\mathrm{d}x}\left(k \dfrac{\mathrm{d}T}{\mathrm{d}x}\right) = 0, \label{eq:heatequationkirchoff}
\end{equation}
where $k = k(T)$ is a temperature-dependent thermal conductivity. Equation (\ref{eq:heatequationkirchoff}) is difficult to solve analytically. If we introduce a new variable by 
\begin{equation}
    \theta(T) = T_0 + \dfrac{1}{k_0} \int_{T_0}^{T} k(s) \:\mathrm{d}s,
\end{equation}
where $T_0$ is a reference temperature (often $T_0 = 0$), $k_0 = k(T_0) $, and $\theta$ is the ``apparent'' temperature, the nonlinear heat conduction equation is transformed into the linear Laplace's equation by
\begin{equation}
    \dfrac{\mathrm{d}^2 \theta}{\mathrm{d}x^2} = 0.
\end{equation}
If the boundary conditions on $T$ can be transformed into linear boundary conditions on $\theta$, the linearized problem can be solved analytically. Then, following an algebraic inverse Kirchoff transform, the exact nonlinear temperature can be reproduced \cite{BagnallMuzychka2014}. For more detail on the Kirchoff transformation, the reader is referred to Refs.~\cite{BagnallMuzychka2014,Joyce1975,BonaniGhione1995}. Here, similar techniques are borrowed for a parabolic PDE of a certain kind. Following the notation of Ref.~\cite{Heinkenschloss1996}, if a new function, $\beta(t)$, is defined by 
\begin{equation}
    \beta(t) = \int_0^t \kappa(s)\:\mathrm{d}s, \label{eq:beta}
\end{equation}
a finite-dimensional weak form can be constructed, written as
\begin{multline}
\int_{\Omega} \dfrac{\partial y^{(j)}_h}{\partial r}v^{(j)}_h\:\mathrm{d}x  + \psi^{(j)}\int_{\Omega}\dfrac{\partial \beta^{(j)}_h}{\partial x}v^{(j)}_h\:\mathrm{d}x \\ =   a \psi^{(j)} \left(u_2^{(j)}v_h^{(j)}(1) - u_1^{(j)}v_h^{(j)}(0)
- \int_{\Omega} \frac{\partial y_h^{(j)}}{\partial x} \frac{\partial v_h^{(j)}}{\partial x}\:\mathrm{d}x \right).
\label{eq:weakform3}
\end{multline}
where $\beta(y^{(j)}) \approx \beta_h^{(j)}$ and $\beta_h^{(j)}$ is defined explicitly by 
\begin{equation}
    \beta(y^{(j)}) \approx \beta_h^{(j)} = \sum_{k = 1}^{N_x} \phi_k \beta^{(j)}_k = \sum_{k = 1}^{N_x} \phi_k\int_0^{Y^{(j)}_k} \kappa (s) \:\mathrm{d}s, \label{eq:betaapprox}
\end{equation}
where $\kappa(Y^{(j)}_k)$ is defined by $\kappa(Y^{(j)}_k) \triangleq \left.\kappa(y)\right|_{y = Y^{(j)}_k}$. This implies that the transformation is applied directly to the coefficients of $y_h$ rather than the entire finite-dimensional function. This creates an approximation, as the degree of the basis functions approximating the nonlinear term lag behind the ``true'' degree. That is, the degree of the product of $\kappa(y_h^{(j)})(\partial y_h^{(j)}/\partial x)$ is larger than that of $\partial \beta_h^{(j)}/\partial x$. 

The key issue is whether the lagged polynomial degree has practical significance. As the solution is being approximated in a polynomial basis, the method inherently relies on the assumption that the true solution is well represented by such a basis. If this assumption fails, the scheme cannot capture the underlying behavior regardless of the degree of the nonlinearity. When the solution is in fact well approximated by a polynomial, it becomes essential that the mesh is sufficiently refined or that the basis functions possess derivatives of adequate order to resolve the characteristic features of the nonlinear term. If the mesh is sufficiently refined, the specific polynomial degree used to approximate the nonlinear term becomes secondary provided the discretization of the domain affords adequate resolution to capture the relevant dynamics.

Note that the ``Kirchoff-like'' transform differs from the original Kirchoff transform. Though similar in nature, in this case, the problem is solved in the primary state variable, which differs from the original Kirchoff transform, where the problem is solved entirely in a transformed variable, and the physical solution is back-calculated via an inverse Kirchoff transform. Ultimately, the transformation presented here serves as an approximation to the original PDE by neglecting higher-order cross-term products of the basis functions, specifically in the nonlinearity. The ``Kirchoff-like'' transform is a way to improve the computational efficiency in generating a solution to the optimal control problem through a linearized approximation, which can be viewed as an analogy to techniques such as mass lumping \cite{ZienkiewiczTaylor2005, Hughes2003, HintonRock1976} in finite element analysis. The weak form given by Eq.~\eqref{eq:weakform3} has the advantage of being linear in a larger set of unknowns, in contrast to the nonlinear dependence on $y$ in Eq.~\eqref{eq:weakform1}, which mitigates the computational complexity of evaluating, and taking derivatives of, the NLP constraints function. If the leading term in the nonlinearity is nonintegrable, other techniques may be required to handle the nonlinearity present.

\subsection{Construction of the Semi-Discrete Form}
\begin{figure}[h!]
      \centering
      \includegraphics[scale=0.9]{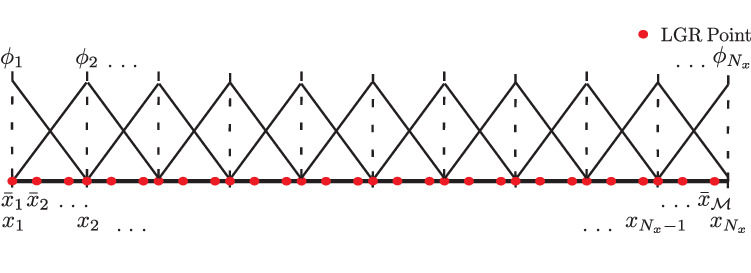}
      \vspace{-1em}
      \caption{A depiction of the spatial discretization with piecewise linear elements and $n_x = 3$ LGR points between each node.}
      \label{fig:linearelements}
   \end{figure}

To construct the semi-discrete form, the discrete functional approximations, $y_h$ and $v_h$, are substituted into the variational form. Expanding $y_h$ and $v_h$ as linear combinations of basis functions and leveraging linearity
leads to the semi-discrete form in terms of $Y_k$ and $\beta_k$, the coefficients of those expansions, as given by
\begin{equation}
            M\mathbf{Y}^{(j)}_r+ \psi^{(j)} N\boldsymbol{\beta}^{(j)} = -a  \psi ^{(j)} A\mathbf{Y}^{(j)} 
        + \psi^{(j)}(\mathbf{e}_{N_x}u_2^{(j)} - \mathbf{e}_{1}u_1^{(j)}), \label{eq:semidiscrete}
\end{equation}
where $\boldsymbol{\beta}^{(j)}= [\beta^{(j)}_1, \ldots,\beta_{N_x}^{(j)}]^T \in \mathbb{R}^{N_x \times 1}$, $\mathbf{Y}^{(j)} = [Y^{(j)}_1, \ldots,Y^{(j)}_{N_x}]^T\in \mathbb{R}^{N_x \times 1}$, and the vector $\mathbf{Y}_r^{(j)} \in \mathbb{R}^{N_x \times 1}$  represents the derivative of the time-dependent coefficients, $\mathbf{Y}^{(j)}$, with respect to the transformed temporal variable. The matrices $M$, $N$, and $A$ can be computed by 
\begin{align}
    M_{ik} &= \int_{\Omega} \phi_k\phi_i\:\mathrm{d}x \in \mathbb{R}^{N_x \times N_x}, \label{eq:M} \\
    N_{ik} &= \int_{\Omega} \frac{\partial \phi_k}{\partial x} \phi_i\:\mathrm{d}x \in \mathbb{R}^{N_x \times N_x}, \label{eq:N} \\
    A_{ik} &= \int_{\Omega} \frac{\partial \phi_k}{\partial x} \frac{\partial \phi_i}{\partial x}\:\mathrm{d}x \in \mathbb{R}^{N_x \times N_x}, \label{eq:A}
\end{align}
and $\mathbf{e}_1$ and $\mathbf{e}_{N_x}$ represent the elementary vectors defined by 
\begin{align}
\mathbf{e}_{1} & = [1,\ldots,0]^T \in \mathbb{R}^{N_x \times 1}, \\
\mathbf{e}_{N_x} & = [0,\ldots,1]^T \in \mathbb{R}^{N_x\times 1}.
\end{align}
It is noted that the integrals in Eqs.~(\ref{eq:M}-\ref{eq:A}) can be computed analytically in certain cases. However, an approximation of the integrals by numerical quadrature avoids this analytical derivation and generalizes to basis functions of higher-degree. The integrals in Eqs.~(\ref{eq:M}-\ref{eq:A}) can be approximated numerically via a Gaussian quadrature by  
\begin{align}
    M_{ik} &= \frac{1}{2}\sum_{p = 1}^{\mathcal{M}} w_p h_p \phi_k(\bar{x}_p)\phi_i(\bar{x}_p), \label{eq:Mquad} \\
    N_{ik} &= \frac{1}{2}\sum_{p = 1}^{\mathcal{M}} w_p h_p \frac{ \partial \phi_k(\bar{x}_p)}{\partial x}\phi_i(\bar{x}_p), \label{eq:Nquad} \\
    A_{ik} &= \frac{1}{2}\sum_{p = 1}^{\mathcal{M}} w_p h_p \frac{ \partial \phi_k(\bar{x}_p)}{\partial x}\frac{ \partial \phi_i(\bar{x}_p)}{\partial x}, \label{eq:Aquad} 
\end{align}
where $\bar{x}_p \in \Omega,\:p = \{1,\ldots,\mathcal{M}\}$, represents the location of a specified set of Gaussian quadrature points on the spatial domain. In this work, a Legendre-Gauss-Radau (LGR) quadrature is utilized with a consistent number of points between each node. LGR quadrature is exact for polynomials of degree $2n_x-2$, where $n_x$ is the user-specified number of LGR points between each node. Thus, $\mathcal{M}$ is the total number of LGR points utilized to approximate the spatial integral and can be computed by $\mathcal{M} = n_x \times (N_x-1)$. The LGR quadrature weight that corresponds to the LGR point at $\bar{x}_p$ is denoted by $w_p$, and $h_p$ indicates the distance between the left and right closest node at $\bar{x}_p$. A depiction of the spatial discretization with piecewise linear elements is provided in Fig.~\ref{fig:linearelements} for $n_x = 3$ LGR points between each node. The PDE is now removed of explicit spatial dependence, and an application of a temporal discretization is required to fully discretize the problem. 

\section{Multi-Interval Flipped LGR Temporal Collocation}\label{sec:fLGR}
\begin{figure}[h]
      \centering
      \includegraphics[scale=0.6]{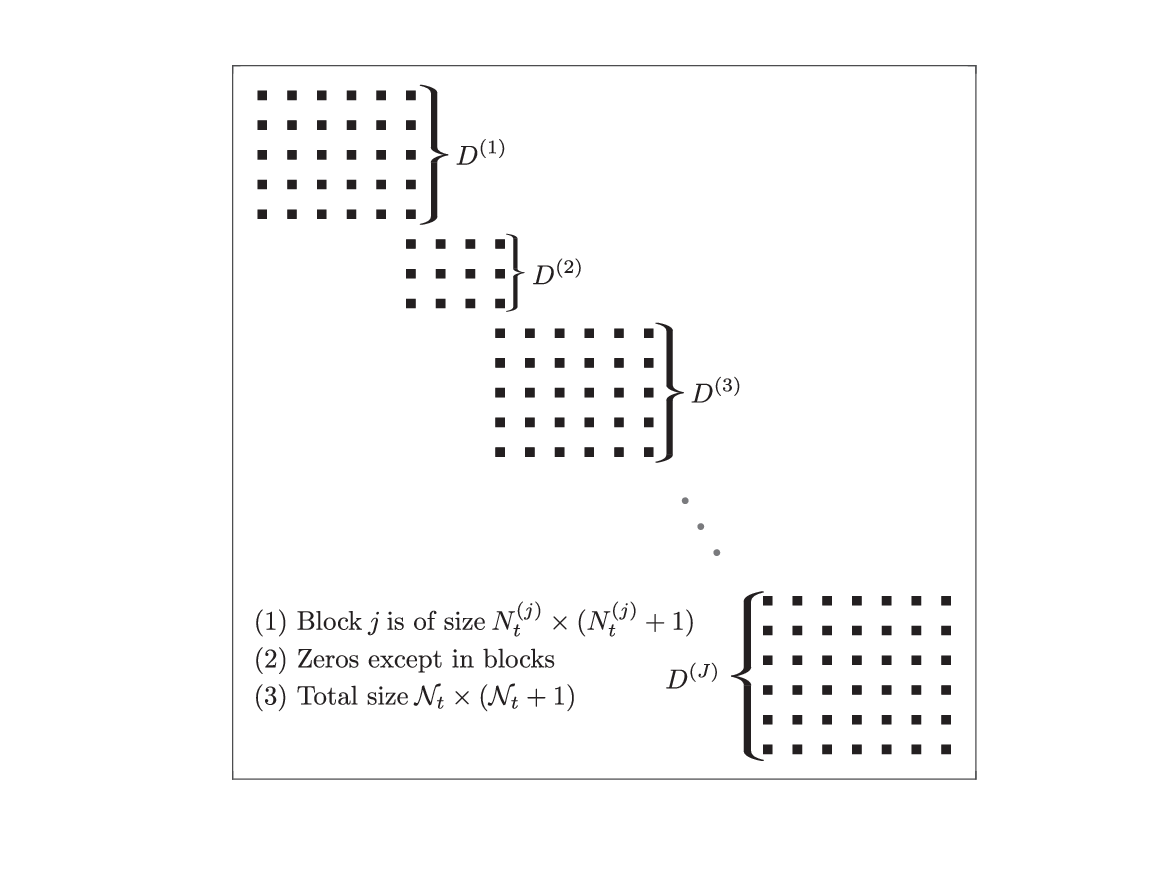}
      \vspace{-1em}
      \caption{Structure of the fLGR temporal differentiation matrix for a mesh with $J$ time intervals.}
      \label{fig:samplediffmat}
   \end{figure}

To discretize in time, we assume that the state may be approximated by a basis of piecewise-continuous Lagrange polynomials supported at the fLGR points and a noncollocated initial point. The fLGR points are the negative roots of the Legendre polynomial $P_{N_t^{(j)}-1}(r) + P_{N_t^{(j)}}(r)$, where $N_t^{(j)}$ is the number of temporal collocation points in interval $j = \{1,\ldots,J\}$. The points exist on the half-open interval $(-1,1]$ and can be utilized to exactly integrate a polynomial $p(r)$ of degree $2N_t^{(j)}-2$ by 
\begin{equation}
    \int_{-1}^{1} p(r)\:\mathrm{d}r = \sum_{i=1}^{N_t^{(j)}} w^{(j)}_i p(r_i),
\end{equation}
where $w_i^{(j)}$, $j=\{1,\ldots J\}$, are the flipped LGR quadrature weights. The spatially discrete, continuous time vector $\mathbf{Y}^{(j)}$ is now approximated in time by a finite-dimensional function $\mathbf{\bar{Y}}^{(j)} \in \mathbb{R}^{N_x \times 1}$, where the finite-dimensional function is used to approximate the behavior of the spatial coefficients, $\mathbf{Y}^{(j)} = [Y^{(j)}_1, \ldots,Y^{(j)}_{N_x}]^T \in \mathbb{R}^{N_x \times 1}$, as they evolve in time. This is done by  
\begin{equation}
    \mathbf{Y}^{(j)} \approx \bar{\mathbf{Y}}^{(j)} = \sum_{n = 0}^{N_t^{(j)}} \bar{\mathbf{Y}}_n^{(j)}L_n^{(j)}(r). \label{eq:stateapprox}
\end{equation}
For notation, a particular index $n$ of $\mathbf{\bar{Y}}_n^{(j)} \in \mathbb{R}^{N_x \times 1}$ refers to the approximation of the vector of spatial coefficients, $\mathbf{Y}$, at the $n^{\mathrm{th}}$ temporal point in interval $j = \{1,\ldots, J\}$, where the Lagrange polynomials, $L^{(j)}_n(r)$, are defined by 
\begin{equation}
     L^{(j)}_n(r) = \prod_{\substack{i=0 \\ i \neq n}}^{N^{(j)}_t} \frac{r - r_{i}^{(j)}}{r_{n}^{(j)}-r_{i}^{(j)}}. \label{eq:Lpoly}
\end{equation}
Though it is noted that other interpolating polynomials may be used (e.g.~Chebyshev, Legendre, Jacobi, etc.), we select Lagrange interpolating polynomials. The use of Lagrange interpolating polynomials in association with the fLGR points can be shown to eliminate the Runge phenomenon \cite{Huntington2007} as well as provide exponential convergence for smooth and well-behaved problems \cite{HagerHou2019}. It is noted that the approximation is identically applied to $\boldsymbol{\beta}^{(j)}$ to yield.
\begin{equation}
    \boldsymbol{\beta}^{(j)} \approx \bar{\boldsymbol{\beta}}^{(j)} = \sum_{n = 0}^{N_t^{(j)}} \bar{\boldsymbol{\beta}}_n^{(j)}L_n^{(j)}(r). \label{eq:stateapprox2}
\end{equation} 
The state approximation in Eq.~\eqref{eq:stateapprox} can be differentiated with respect to the transformed temporal variable and applied at the collocation points, $r_i$, $i=\{1,\ldots,N_t^{(j)}\}$, by 
\begin{equation}
    \dfrac{\mathrm{d} \mathbf{Y}^{(j)}(r_i)}{\mathrm{d} r} \approx  \frac{\mathrm{d} \bar{\mathbf{Y}}^{(j)}(r_i)}{\mathrm{d} r} = \sum_{n=0}^{N_t^{(j)}}  \bar{\mathbf{Y}}_n^{(j)} \frac{\mathrm{d} L^{(j)}_n(r_i)}{\mathrm{d} r}. \label{eq:derivapprox1}
\end{equation}
An equivalent expression for the derivative of the state approximation in Eq.~\eqref{eq:derivapprox1} is provided as 
\begin{equation}
   \frac{\mathrm{d} \bar{\mathbf{Y}}^{(j)}(r_i)}{\mathrm{d} r} = \sum_{n=0}^{N_t^{(j)}}  \bar{\mathbf{Y}}_n^{(j)} \frac{\mathrm{d} L^{(j)}_n(r_i)}{\mathrm{d} r} = \sum_{n=0}^{N_t^{(j)}} D_{in}^{(j)}\bar{\mathbf{Y}}_n^{(j)},
\end{equation}
where $D_{in}^{(j)}$ is the fLGR differentiation matrix in interval $j = \{1,\ldots, J\}$, defined as  
\begin{equation}
    D_{in}^{(j)} =  \left[\frac{\mathrm{d}L_n^{(j)}(r)}{\mathrm{d} r}\right]_{r_i^{(j)}} \quad \in \mathbb{R}^{N_t^{(j)} \times (N_t^{(j)} +1) }. 
\end{equation}
If we define the matrix
\begin{equation}
    \bar{Y}^{(j)} = \begin{bmatrix}
        \bar{\mathbf{Y}}_{0}^{(j)} & \bar{\mathbf{Y}}_{1}^{(j)} & \ldots & \bar{\mathbf{Y}}_{N_t^{(j)}}^{(j)} \\
    \end{bmatrix} \quad \in \mathbb{R}^{N_x \times (N_t^{(j)}+1)}, \label{eq:intstatematrix}
\end{equation}
through utilizing the same variable in the NLP at temporal mesh interval intersections, a state matrix can be constructed by 
\begin{equation}
    \bar{Y} = \begin{bmatrix}
        \bar{Y}^{(1)}& \bar{Y}^{(2)}_{1:N^{(2)}_t} & \ldots & \bar{Y}^{(J)}_{1:N^{(J)}_t}
    \end{bmatrix}^T \quad \in \mathbb{R}^{(\mathcal{N}_t+1) \times N_x},
\end{equation}
where $\mathcal{N}_t$ is the total number of collocation points in the temporal dimension given by 
\begin{equation}
    \mathcal{N}_t = \sum_{j=1}^J N_t^{(j)}.
\end{equation}
It is noted that the columns of $\bar{Y}$ are denoted by $\mathbf{\bar{Y}}_1,$ $\mathbf{\bar{Y}}_2, \ldots, \mathbf{\bar{Y}}_{N_x} \in \mathbb{R}^{(\mathcal{N}_t +1) \times 1}$. The construction of the state matrix, $\bar{Y}$, allows for the construction of a matrix product that formulates a fully-discrete form of Eq.~\eqref{eq:semidiscrete}. The fully-discrete form can be constructed by 
\begin{equation}
        M\left(D_t \bar{Y}\right)^T+\boldsymbol{\alpha} \odot (N\bar{\beta}^T) = -a \boldsymbol{\alpha} \odot  (A(\bar{Y}_{1:\mathcal{N}_t,:})^T) +\boldsymbol{\alpha} \odot (\mathbf{e}_{N_x}\mathbf{U}_2 - \mathbf{e}_{1}\mathbf{U}_1),  \label{eq:fullydiscrete} 
\end{equation}
where the operator $\mathbf{f} \odot G$ represents the element-by-element product of vector $\mathbf{f}$ with matrix $G$. The matrix, $\bar{\beta} \in \mathbb{R}^{\mathcal{N}_t \times N_x}$, can be expressed as 
\begin{equation}
      \bar{\beta} = \begin{bmatrix}
        \bar{\beta}^{(1)}_{1:N^{(1)}_t} & \bar{\beta}^{(2)}_{1:N^{(2)}_t} & \ldots & \bar{\beta}^{(J)}_{1:N^{(J)}_t}
    \end{bmatrix}^T, \label{eq:discretebeta}
\end{equation}
where $\bar{\beta}^{(j)}$ is defined similarly to Eq.~(\ref{eq:intstatematrix}). Additionally, the vector of temporal interval mappings evaluated at the collocation points, $\boldsymbol{\alpha}$, is given by  
\begin{equation}
    \boldsymbol{\alpha} = \begin{bmatrix}
        \boldsymbol{\psi}^{(1)} &  \boldsymbol{\psi}^{(2)} & \ldots &  \boldsymbol{\psi}^{(J)}
    \end{bmatrix} \quad \in \mathbb{R}^{1 \times \mathcal{N}_t}, 
    \end{equation}
 where $ \boldsymbol{\psi}^{(j)}$ is defined by 
    \begin{equation}
         \boldsymbol{\psi}^{(j)} = \begin{bmatrix}
             \psi^{(j)} & \psi^{(j)} & \ldots & \psi^{(j)}
         \end{bmatrix} \quad \in \mathbb{R}^{1 \times N_t^{(j)}}.
    \end{equation}
    In other words, $\psi^{(j)}$ appears $N_t^{(j)}$ times in $\boldsymbol{\alpha}$ corresponding to each collocation point in interval $j$. The control is defined as a set of discrete variables located at the collocation points.
    \begin{align}
    u_1^{(j)} &\approx \mathbf{U}_1^{(j)} = \begin{bmatrix}
        U_{1,1}^{(j)} & U_{1,2}^{(j)} & \ldots & U_{1,N_t^{(j)}}^{(j)} 
        \end{bmatrix}\quad \in \mathbb{R}^{1 \times N_t^{(j)}}, \\
        u_2^{(j)} &\approx \mathbf{U}_2^{(j)} = \begin{bmatrix}
        U_{2,1}^{(j)} & U_{2,2}^{(j)} & \ldots & U_{2,N_t^{(j)}}^{(j)}     
    \end{bmatrix}\quad \in \mathbb{R}^{1 \times N_t^{(j)}},
\end{align}
where the notation $U^{(j)}_{k,i}$ indicates the $k^\mathrm{th}$ control variable located at $r^{(j)}_i$. Control vectors $\mathbf{U}_1$ and $\mathbf{U}_2$ can then be constructed by 
\begin{align}
    \mathbf{U}_1 = \begin{bmatrix}
        \mathbf{U}_1^{(1)} & \mathbf{U}_1^{(2)} & \ldots & \mathbf{U}_1^{(J)}
    \end{bmatrix} \quad \in \mathbb{R}^{1 \times \mathcal{N}_t}, \\
        \mathbf{U}_2 = \begin{bmatrix}
        \mathbf{U}_2^{(1)} & \mathbf{U}_2^{(2)} & \ldots & \mathbf{U}_2^{(J)}
    \end{bmatrix} \quad \in \mathbb{R}^{1 \times \mathcal{N}_t}.
\end{align}
The temporal differentiation matrix, $D_t$, is constructed by 
\begin{equation}
D_t = 
    \begin{bmatrix}
      D^{(1)} & 0 & \ldots & 0 \\
      0 & D^{(2)} & \ddots & \vdots \\
      \vdots & \ddots & \ddots & \vdots \\
      0& \ldots & \ldots & D^{(J)}
    \end{bmatrix} \quad \in \mathbb{R}^{\mathcal{N}_t \times (\mathcal{N}_t+1)}, \label{eq:diffmat}
\end{equation}
where it is noted that there exists one overlapping column over each interval differentiation matrix due to enforcement of continuity in the state through the use of the same variable in the NLP: 
\begin{equation}
    \bar{\mathbf{Y}}^{(j)}_{N_t^{(j)}} = \bar{\mathbf{Y}}^{(j+1)}_{0}, \quad j = \{1,\ldots,J\}.
\end{equation}
An example temporal differentiation matrix $D_t$ is provided in Fig.~\ref{fig:samplediffmat} for clarity.
The discretization of the initial conditions can be completed through the enforcement of the initial 
data at the points that correspond to the initial time. Explicitly, this is done by 
\begin{equation}
    \bar{Y}_{0k} = q(x_k,t_0), \quad k =\{1,\ldots,N_x\}.  \label{eq:InitialConditions}
\end{equation}
Thus, the only remaining procedure to completely discretize the optimal control problem is the discretization of the objective. This can be done by 
\begin{equation}
\mathcal{J} \approx \sum_{i=1}^{\mathcal{N}_t}\omega_{i}\left[ \frac{1}{2}\sum_{p= 1}^{\mathcal{M}} h_p w_p \mathcal{L}\left(\bar{x}_p,t_i,\left(\sum_{k=1}^{N_x} \phi_{k}(\bar{x}_p)\bar{Y}_{ik}\right)\right) + \mathcal{P}(t_i, U_{1,i}, U_{2,i}) \right],\label{eq:discreteobjective}
\end{equation}
where $\omega_i$ is the $i^{\mathrm{th}}$ component of $\boldsymbol{\omega} \in \mathbb{R}^{1 \times \mathcal{N}_t} $ where $\boldsymbol{\omega}$ is a vector containing the corresponding quadrature weights on $r^{(j)}_1, r^{(j)}_2, \ldots, r^{(j)}_{N_t^{(j)}}\:\forall\:j = \{1,\ldots, J\}$, or, more explicitly:
\begin{equation}
    \boldsymbol{\omega} = \begin{bmatrix}
        \mathbf{w}_t^{(1)} & \mathbf{w}_t^{(2)} & \ldots & \mathbf{w}_t^{(J)}
    \end{bmatrix}\quad \in \mathbb{R}^{1 \times \mathcal{N}_t},
\end{equation}
where 
\begin{equation}
    \mathbf{w}_t^{(j)} = \begin{bmatrix}
        w_1^{(j)} & w_2^{(j)} & \ldots & w_{N_t^{(j)}}^{(j)}
    \end{bmatrix} \quad \in \mathbb{R}^{1 \times N_t^{(j)}}.
\end{equation}
For completeness, $t_i$ represents the time at the corresponding collocation point.

\begin{figure}[h]
    \subfloat[fLGR discretization\label{fig:fLGRBC}]{\includegraphics[scale=0.4]{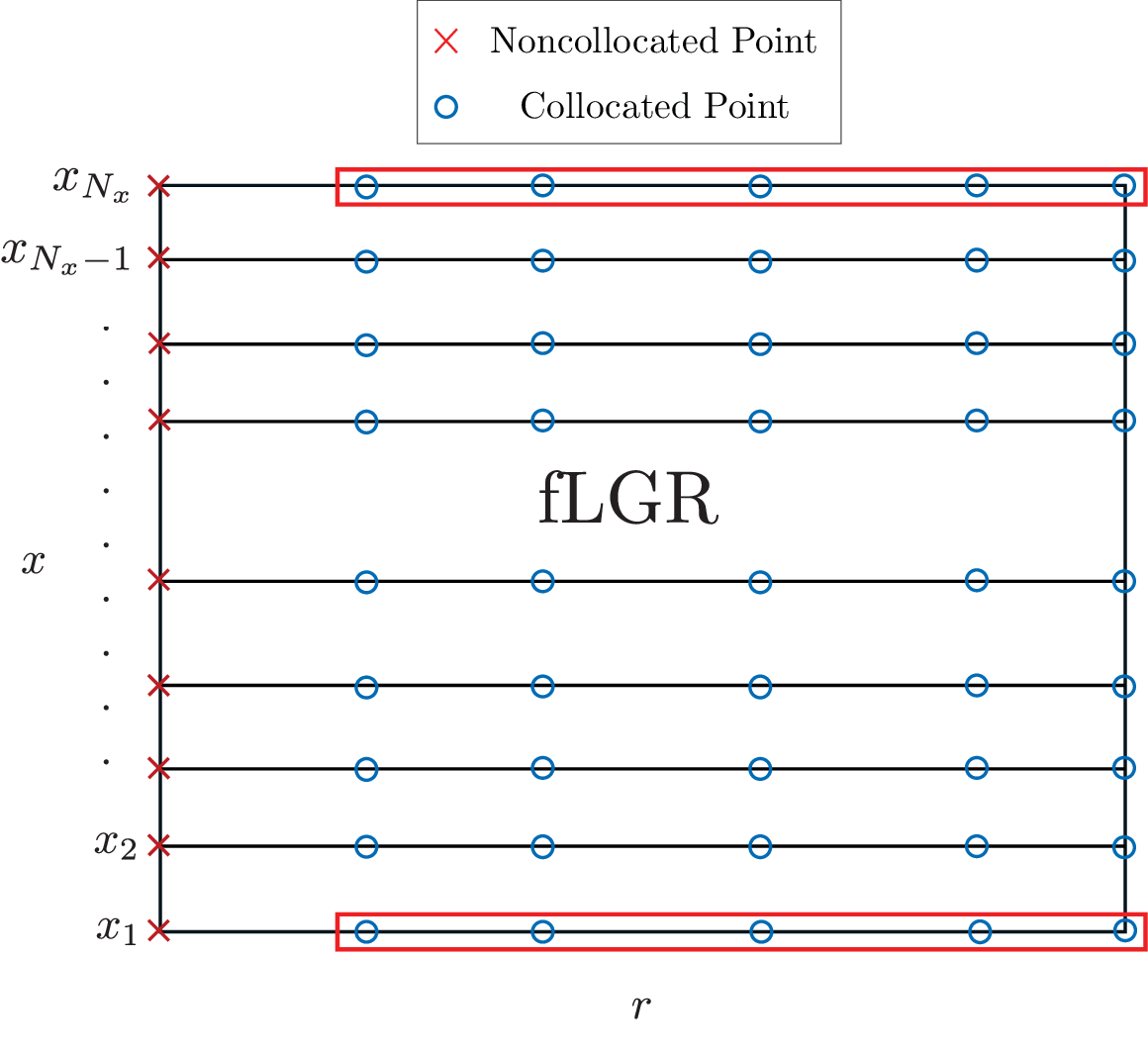}}
    \hfill
    \subfloat[Standard LGR discretization\label{fig:LGRBC}]{\includegraphics[scale=0.4]{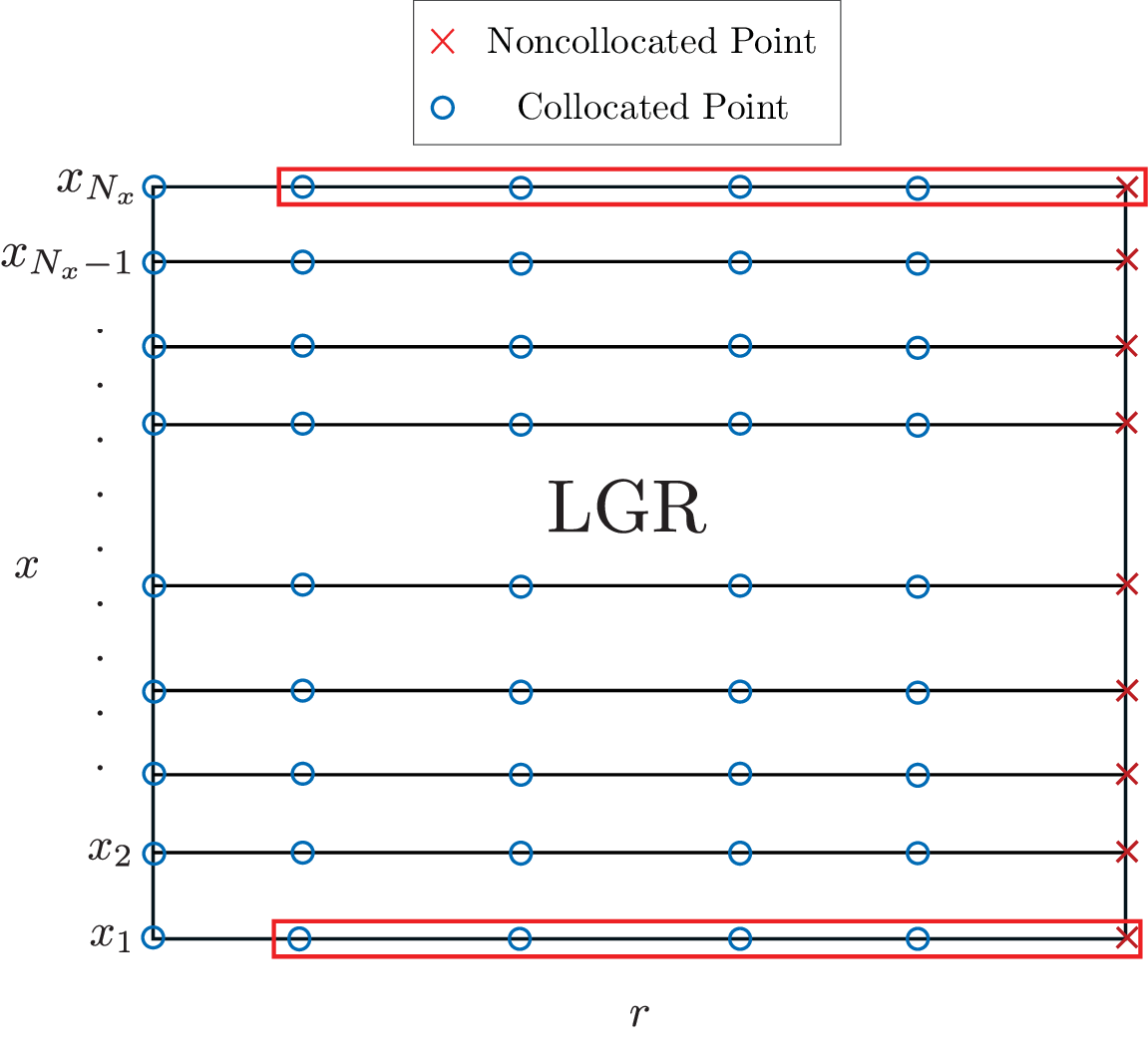}}
    \caption{A comparison between the fLGR and LGR temporal discretizations and the locations of the applied boundary conditions (boxed). \label{fig:LGRComp}}
\end{figure}

Generally, in optimal control problems governed by ordinary differential equations, the standard LGR and fLGR schemes may be used interchangeably to achieve similar solutions \cite{GargPatterson2010}. In the context of optimal control problems governed by PDEs, however, the use of the fLGR points as opposed to the standard LGR points allows for the satisfaction of the boundary conditions at the terminal time. Because the fLGR points possess a collocated point at the terminal time, the boundary conditions are related to the dynamics through the discretized form in Eq.~\eqref{eq:fullydiscrete}. If the standard LGR points are utilized, the boundary condition cannot be satisfied through the discretized form, as the terminal time exists as a noncollocated point. A visual representation of the comparison between the two meshes presented by the schemes and the location of the points at which the boundary conditions are applied is provided in Fig.~\ref{fig:LGRComp}. It is noted that the time discretizations in Fig.~\ref{fig:LGRComp} are shown in the global (one interval) sense for visual simplicity. 
It is evident in Fig.~\ref{fig:LGRBC} that the boundary condition at the terminal time cannot be related to the interior dynamics due to the noncollocated point at the terminal time. The initial time is not boxed as it is constrained by the initial condition. Thus, in a sense, the fLGR discretization is constructed in a manner that is ideal for time-dependent PDE constrained optimization problems, as the initial condition is typically a fixed parameter of a system, and the general goal of most time-dependent PDE optimal control problems is to drive the state from a specified initial condition to a desired profile through the control. It is possible that other alternatives exist to address this phenomenon. One may consider, for example, a similar approach to the one taken in Ref.~\cite{DaviesDennis2024}, where collocation points are added to the scheme to satisfy boundary conditions at both boundaries. However, it is anticipated that this technique will lead to similar downfalls as those presented in Ref.~\cite{DaviesDennis2024}, where the NLP search space was negatively affected and avenues for user error were introduced. Thus, the fLGR scheme is a preferable scheme in its application to optimal control problems governed by PDEs, as the location of the collocation points are well-suited to satisfy dynamics/boundary relations. 

\section{Transcription of the Optimal Control Problem into an NLP}\label{sec:NLP}
The NLP associated with the discretized optimal control problem is provided as follows. Determine the decision vector, $\mathbf{z}$, that minimizes 
\begin{equation}
    \Phi(\mathbf{z}),
\end{equation}
subject to the constraints 
\begin{equation}
   \mathbf{F}_{\text{min}}  \leq \mathbf{F}(\mathbf{z}) \leq \mathbf{F}_{\text{max}},
\end{equation}
and the variable bounds 
\begin{equation}
    \mathbf{z}_{\text{min}} \leq \mathbf{z} \leq \mathbf{z}_{\text{max}}.
\end{equation}
The size of the NLP is determined by the number of temporal collocation points, mesh intervals, and spatial nodes utilized to complete the discretization. Due to the added dimensionality of PDEs, the size of the decision vector scales by a product of the number of temporal and spatial points. Though the size of the decision vector can change, the structure of the NLP is invariant of the mesh size. Additionally, as the size of the mesh grows, the resulting sparsity of the NLP grows. The decision vector can be organized by 
\begin{equation}
\mathbf{z} = 
\begin{bmatrix}
    \mathbf{\bar{Y}}(:)^T &
    \mathbf{U}_1 &
    \mathbf{U}_2 &
    t_0 &
    t_f
    \end{bmatrix} \quad \in \mathbb{R}^{1 \times \left[(\mathcal{N}_t+1)N_x + 2\mathcal{N}_t + 2\right]},
\end{equation}
where the operator $\mathbf{\bar{Y}}(:)$ represents the reorganization of matrix $\bar{Y}$ into a single column vector. More explicitly, this is given by 
\begin{equation}
    {\mathbf{\bar{Y}}}(:) = \begin{bmatrix}
        \mathbf{\bar{Y}}_{1}^T & \mathbf{\bar{Y}}_{2}^T & \ldots & \mathbf{\bar{Y}}_{N_x}^T
    \end{bmatrix}^T, \quad \in \mathbb{R}^{(\mathcal{N}_t+1)N_x \times 1}. \label{eq:Ycolon}
\end{equation}
The bounds on the constraint functions $\mathbf{F}(\mathbf{z})$ given by $\mathbf{F}_{\text{min}}$ and $\mathbf{F}_{\text{max}}$ along with the variable bounds $\mathbf{z}_{\text{min}}$ and $\mathbf{z}_{\text{max}}$ are specified by the user. The NLP objective is computed, $\Phi(\mathbf{z}) = \Phi$, by applying Eq.~(\ref{eq:objective}) at the discrete collocation points (Eq.~(\ref{eq:discreteobjective})) where 
\begin{equation}
    \mathcal{J} \approx \Phi = \sum_{i=1}^{\mathcal{N}_t}\omega_{i}\left[ \frac{1}{2}\sum_{p= 1}^{\mathcal{M}} h_p w_p \mathcal{L}\left(\bar{x}_p,t_i,\left(\sum_{k=1}^{N_x} \phi_{k}(\bar{x}_p)\bar{Y}_{ik}\right)\right) + \mathcal{P}(t_i, U_{1,i}, U_{2,i}) \right]. \label{eq:Phi}
\end{equation}
The constraint functions are assembled by 
\begin{equation}
    \mathbf{F} = \begin{bmatrix}
        \bar{\mathbf{Y}}_{0}^{(1)}-q(\mathbf{x},t_0) \\
         \mathrm{Eq.}~\eqref{eq:fullydiscrete}(:)
    \end{bmatrix}\:\: \in \mathbb{R}^{(\mathcal{N}_tN_x + N_x) \times 1},\label{eq:constraintfuncs}
\end{equation}
where $\mathbf{x} = \begin{bmatrix}
    x_1 & x_2 & \ldots& x_{N_x}
\end{bmatrix}^T \in \mathbb{R}^{N_x \times 1}$. 
It is important to notice that Eq.~(\ref{eq:discretebeta}) is not explicitly enforced as an NLP constraint. Rather, $\bar{\beta}$ is computed based on the coefficients of the state matrix upon each NLP iteration. 
\subsection{Derivatives of the NLP Objective and Constraint Functions}
For most NLP solvers, the gradient of the objective function with respect to the decision vector is required. The gradient of the objective function with respect to the decision vector is written as 
\begin{equation}
    \dfrac{\partial \Phi}{\partial \mathbf{z}} = \begin{bmatrix}
        \dfrac{\partial \Phi}{\partial \bar{\mathbf{Y}}(:)^T} & \dfrac{\partial \Phi}{\partial \mathbf{U}_1} & \dfrac{\partial \Phi}{\partial \mathbf{U}_2} &  \dfrac{\partial \Phi}{\partial t_0} & \dfrac{\partial \Phi}{\partial t_f} 
    \end{bmatrix} \quad \in \mathbb{R}^{1 \times \left[(\mathcal{N}_t+1)N_x + 2\mathcal{N}_t + 2\right]}.
\end{equation}
As the objective is typically an integral cost that is a function of the state, computation of the state approximation $y_h$ at the quadrature points is usually required. As the mesh is determined prior to the solution of the NLP, the value of the basis functions (and its derivatives) at the quadrature points can be computed prior to solving the NLP. Thus, the objective may be quickly evaluated, as the quadrature points and quadrature formulation are determined prior to the solution of the NLP and the discretized objective is solely dependent on the NLP variables. The gradient of the objective function may then either be computed analytically from the problem information, or algorithmically via derivative approximation techniques such as the bi-complex step method \cite{AgamawiRao2020} or automatic differentiation. Typically, objectives for PDE-constrained optimal control problems take on the form of tracking-type objectives. Generally, the form of which can be expressed as
\begin{equation}
    \mathcal{J} = \dfrac{1}{2}\int_{t_0}^{t_f}\int_0^1 (y(x,t)-y_d(x,t))^2\:\mathrm{d}x\mathrm{d}t + \frac{1}{2}\int_{t_0}^{t_f} || \mathbf{u}(t) ||^2\:\mathrm{d}t, \label{eq:trackingobj}
\end{equation}
where $y_d(x,t)$ is some desired function to be tracked and $\mathbf{u}(t)$, in the case of boundary control, can either be of size $u(t) \in \mathbb{R}$ for one control input or $\mathbf{u}(t) \in \mathbb{R}^2$ for input on both boundaries (i.e.~$\mathbf{u}(t) = [u_1(t), u_2(t)]^T \in \mathbb{R}^2$). Assuming the more general case where control appears on both boundaries, the form of Eq.~\eqref{eq:Phi} becomes 
\begin{equation}
    \Phi = \frac{1}{2}\sum_{i=1}^{\mathcal{N}_t}\omega_{i}\left[ \frac{1}{2}\sum_{p= 1}^{\mathcal{M}} \left(h_p w_p \left[\left(\sum_{k=1}^{N_x} \phi_{k}(\bar{x}_p)\bar{Y}_{ik}\right) - y_d(\bar{x}_p,t_i)\right]^2 \right)+ U_{1,i}^2 +U_{2,i}^2 \right].
\end{equation}
The gradient with respect to the state vector can be computed as follows. A matrix, $\nabla_{\bar{Y}} \Phi \in \mathbb{R}^{(\mathcal{N}_t+1) \times N_x}$, can be constructed that contains the elements of $\partial \Phi/\partial \bar{\mathbf{Y}}(:)^T \in \mathbb{R}^{1 \times (\mathcal{N}_t+1)N_x}$ as
\begin{equation}
    \nabla_{\bar{Y}_{ij}}\Phi =  \frac{1}{2}\omega_{i}\sum_{p= 1}^{\mathcal{M}} h_p w_p \phi_j(\bar{x}_p)\left[\left(\sum_{k=1}^{N_x} \phi_{k}(\bar{x}_p)\bar{Y}_{ik}\right) - y_d(\bar{x}_p,t_i)\right],  
\end{equation}
 for $i = \{1,\ldots,\mathcal{N}_t\},\:j = \{1,\ldots, N_x\}$. The gradient of the objective with respect to the initial state, $y(x,t_0)$, is a row of zeros as the initial time is noncollocated (i.e.~$\nabla_{\bar{Y}_{0,:}} \Phi = \mathbf{0}^T \in \mathbb{R}^{1 \times N_x}$). The portion of the gradient corresponding to the state can then be given by 
\begin{equation}
    \frac{\partial \Phi}{\partial \bar{\mathbf{Y}}(:)^T} = \left(\nabla_{\bar{Y}}\Phi(:)\right)^T. 
\end{equation}
The gradient with respect to the controls $\mathbf{U}_1$ and $\mathbf{U}_2$ are more straightforward and are given by 
\begin{align}
    \frac{\partial \Phi}{\partial \mathbf{U}_1} &= \boldsymbol{\omega} \odot \mathbf{U}_1, \\
    \frac{\partial \Phi}{\partial \mathbf{U}_2} &= \boldsymbol{\omega} \odot \mathbf{U}_2.
\end{align}
The partial derivatives with respect to the initial and final times are scalar-valued, and the only entity that is a function of the initial and final times is $\omega_i$. The partial derivatives are given by 
\begin{align}
    \frac{\partial \Phi}{\partial t_0} &= \frac{1}{2}\sum_{i=1}^{\mathcal{N}_t}\frac{\partial \omega_{i}}{\partial t_0} \left[ \frac{1}{2}\sum_{p= 1}^{\mathcal{M}} \left( h_p w_p \left[\left(\sum_{k=1}^{N_x} \phi_{k}(\bar{x}_p)\bar{Y}_{ik}\right) - y_d(\bar{x}_p,t_i)\right]^2 \right) + U_{1,i}^2 +U_{2,i}^2 \right], \\
    \frac{\partial \Phi}{\partial t_f} &= \frac{1}{2}\sum_{i=1}^{\mathcal{N}_t}\frac{\partial \omega_{i}}{\partial t_f}\left[ \frac{1}{2}\sum_{p= 1}^{\mathcal{M}} \left( h_p w_p \left[\left(\sum_{k=1}^{N_x} \phi_{k}(\bar{x}_p)\bar{Y}_{ik}\right) - y_d(\bar{x}_p,t_i)\right]^2 \right) + U_{1,i}^2 +U_{2,i}^2 \right],
\end{align}
where $\partial \omega_{i}/\partial t_0$ and $\partial \omega_{i}/\partial t_f$ are given by 
\begin{align}
    \frac{\partial \omega_{i}}{\partial t_0} &= -\frac{1}{2} \omega_{i} \frac{\partial \tau}{\partial t}, \\
    \frac{\partial \omega_{i}}{\partial t_f} &= \frac{1}{2} \omega_{i} \frac{\partial \tau}{\partial t}. 
\end{align} 
Thus, the gradient of the NLP objective with respect to the decision vector is complete for objective functionals of the form in Eq.~(\ref{eq:trackingobj}).

\begin{figure}[h]
      \centering
      \includegraphics[scale=0.325]{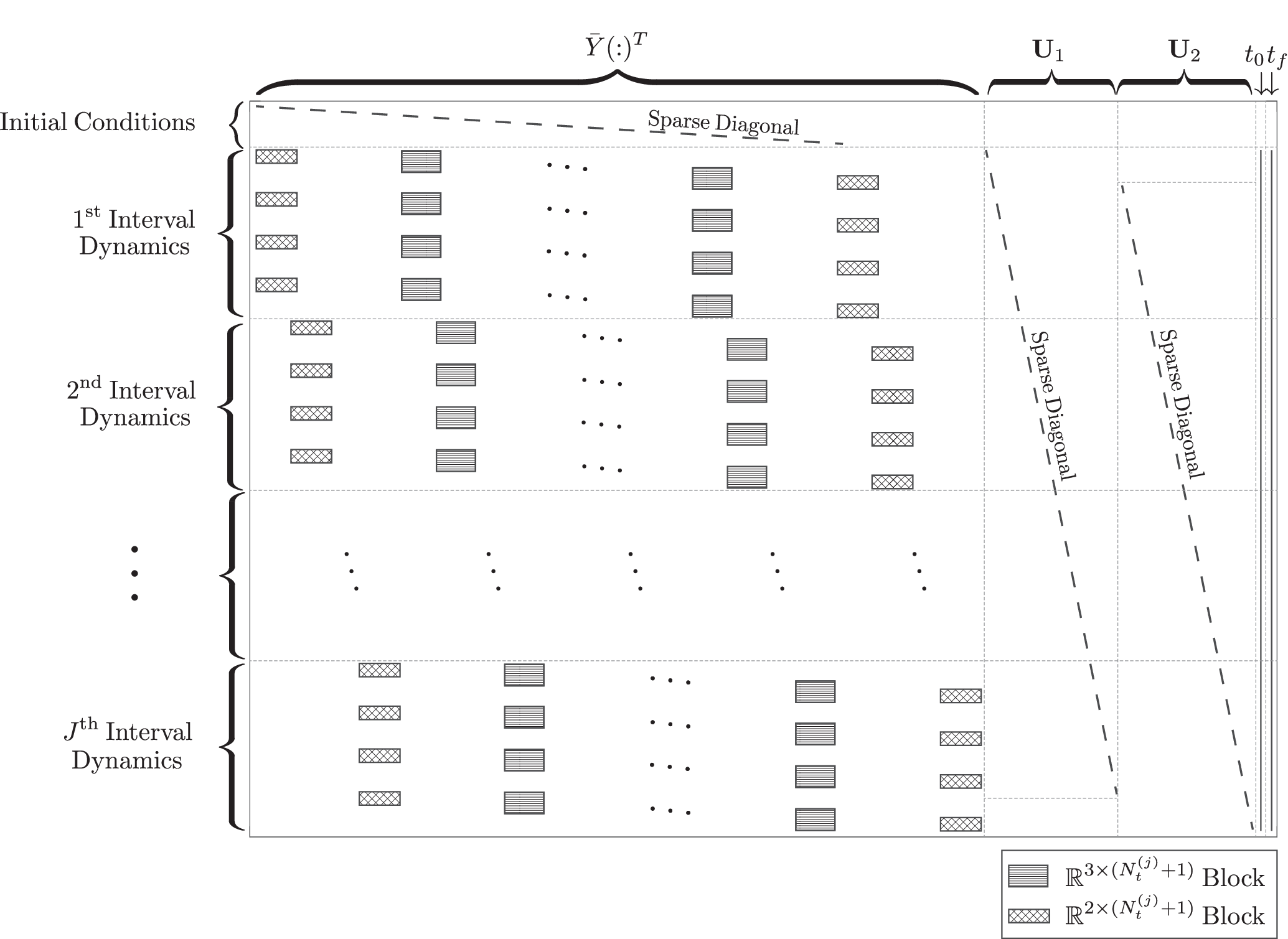}
      \caption{Resulting structure of the NLP constraint Jacobian with P-1 Lagrange elements.}
      \label{fig:JacobianStructure}
   \end{figure}

The Jacobian of the NLP constraints function (for the optimal control problem described in Eqs.~(\ref{eq:objective}-\ref{eq:initialcondition})) is defined by 
\begin{equation}
    \nabla_{\mathbf{z}}\mathbf{F}(\mathbf{z}) = \begin{bmatrix}
        \dfrac{\partial \mathbf{F}}{\partial \bar{\mathbf{Y}}(:)^T} & \dfrac{\partial \mathbf{F}}{\partial \mathbf{U}_1} &  \dfrac{\partial \mathbf{F}}{\partial \mathbf{U}_2} & \dfrac{\partial \mathbf{F}}{\partial t_0} & \dfrac{\partial \mathbf{F}}{\partial t_f}
    \end{bmatrix} \quad \in \mathbb{R}^{(\mathcal{N}_tN_x+N_x) \times \left[(\mathcal{N}_t+1)N_x + 2\mathcal{N}_t + 2\right]}. \label{eq:jacobian}
\end{equation}
The Jacobian of the NLP constraints function (constraint Jacobian) is a sparse matrix dependent upon the mesh. That is, as the number of points in the mesh grows, the number of zero elements in Eq.~(\ref{eq:jacobian}) further exceeds the number of its non-zero elements. The structure of the resulting constraint Jacobian is identical regardless of the mesh size. A graphic depicting the structure of the constraint Jacobian for the optimal control problem in Eqs.~(\ref{eq:objective}-\ref{eq:initialcondition}) is provided in Fig.~\ref{fig:JacobianStructure}.
The diagonal marked by the initial conditions is sparse due to the construction of $\bar{Y}$ and its rearrangement into $\bar{\mathbf{Y}}(:)^T$ (only $N_x$ total non-zero elements). Note, in Eq.~(\ref{eq:Ycolon}), the points corresponding to the initial conditions as highlighted by Eq.~(\ref{eq:InitialConditions}) only appear in increments of $\mathcal{N}_t+1$ (column-wise). Explicitly, the only non-zero elements on the sparse diagonal corresponding to the initial conditions appear at $\nabla_{\mathbf{z}}\mathbf{F}_{i,1+(\mathcal{N}_t+1)(i-1)}$, $\forall\:i = \{1,\ldots,N_x\}$. For the controls $\mathbf{U}_1$ and $\mathbf{U}_2$, the sparsity of their respective diagonals are determined by their sparse appearance in $\mathbf{F}$ (Eq.~(\ref{eq:constraintfuncs})). In other words, the controls only appear in the Jacobian in rows that correspond to $\mathbf{F}$ at the boundaries. It can be shown that the non-zero elements of $\nabla_{\mathbf{U}_1}\mathbf{F}$ and $\nabla_{\mathbf{U}_2}\mathbf{F}$ appear in increments of $N_x$ (row-wise). The non-zero elements appear at $\nabla_{\mathbf{U}_1}\mathbf{F}_{i+N_x(i),i}$, $\forall\:i = \{1,\ldots,\mathcal{N}_t\}$, and $\nabla_{\mathbf{U}_2}\mathbf{F}_{N_x(i)+N_x,i}$, $\forall\:i = \{1,\ldots,\mathcal{N}_t\}$. Additionally, for problems with one control, the columns corresponding to the relevant control are maintained, and the columns corresponding to the absent control are removed.   

The organization of the non-zero elements in the Jacobian corresponding to the partials of the dynamic constraints with respect to the state matrix is a side effect of the organization of $\bar{Y}$ and Eq.~(\ref{eq:fullydiscrete}). The coupling in space and time is demonstrated in Eq.~(\ref{eq:fullydiscrete}), and as a result, some combination of the structure of the differentiation matrices in space and time is expected. The size of each block indicated in the legend of Fig.~\ref{fig:JacobianStructure} highlights this combination. For a particular collocation point on the mesh, the time derivative of the state at that point is evaluated using $N_t^{(j)}+1$ neighboring points within the interval. If the collocation point is located on the interior of the domain, the total number of non-zero locations in a particular row of $M,\:N,\:\mathrm{and}\:A$ is at most three for P-1 Lagrange elements. If the collocation point is located on the boundary of the domain, the total number of non-zero locations in a particular row of $M,\:N,\:\mathrm{and}\:A$ is at most two for P-1 Lagrange elements. The number of times in which the structure is repeated (vertically) in a particular interval corresponds to the number of temporal collocation points in that interval, $N_t^{(j)}$. Lastly, the number of times in which the $3 \times (N_t^{(j)}+1)$ block appears (horizontally within an interval) is given by $N_x-2$, or the number of interior spatial points.

Due to the number of constraints in Eq.~(\ref{eq:constraintfuncs}) and their potential complexity, it is desirable to avoid the computation of the constraint Jacobian by hand. Linearization of the PDE system allows the user to perform numerical integration prior to the optimization procedure, and as a result, the matrices $M$, $N$, $A$ can be computed a priori. Consequently, one can supply derivatives of the NLP constraint function via available automatic differentiation software. In this work, the open source automatic differentiation software \textbf{adigator} \cite{weinstein2017algorithm} was used to supply derivatives of the NLP constraints and objective functions.

\section{Numerical Examples}\label{sec:NumericalExamples}
In this section, two examples are provided to demonstrate the presented method in practice. The first problem analyzes the optimal control of the viscous Burgers' equation with boundary control, and the second problem examines the optimal boundary control of a nonlinear heat equation. All numerical solutions were computed using MATLAB R2023b on a 2023 Apple M2 Ultra Mac Studio with the open source NLP solver {\em IPOPT} \cite{BieglerZavala2008} in full Newton mode with the limited memory Hessian approximation. Convergence analysis was performed using the exact Hessian as computed by \textbf{adigator} for enhanced NLP accuracy at the expense of computational speed. Solutions to both problems are compared to existing solutions in the literature. 

\subsection{Burgers' Equation Tracking Problem}
\begin{figure}[h]
    \centering
    \includegraphics[scale=0.325]{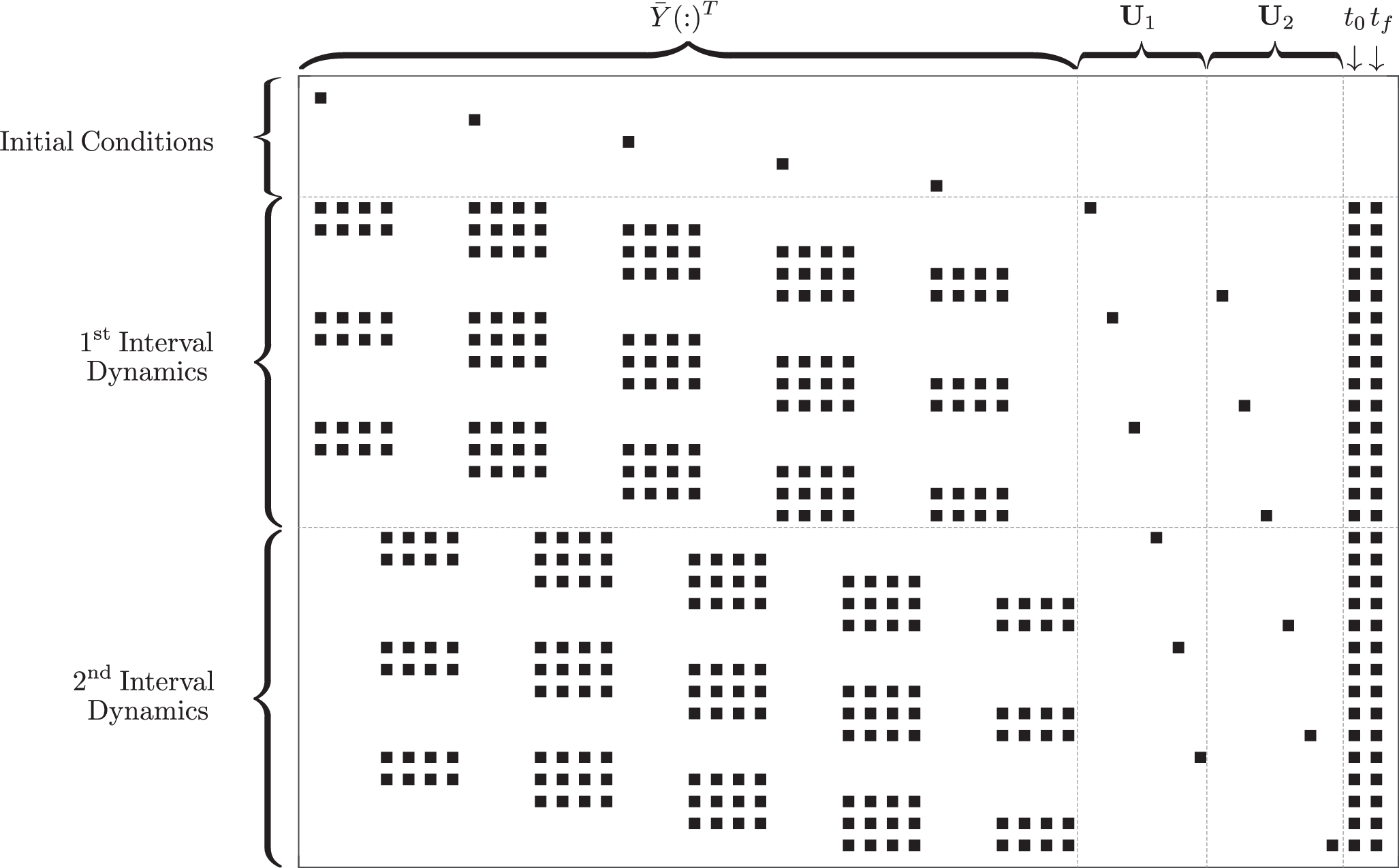}
    \caption{The NLP constraint Jacobian for the Burgers' equation tracking problem with P-1 Lagrange elements and $N_t^{(1)} = N_t^{(2)} = 3$, $J = 2$, $N_x = 5$.}
    \label{fig:constrJacobBurgers}
\end{figure}

\begin{figure}[h!]
    \subfloat[Optimal state \label{fig:burgerstate}]{\includegraphics[scale=0.4]{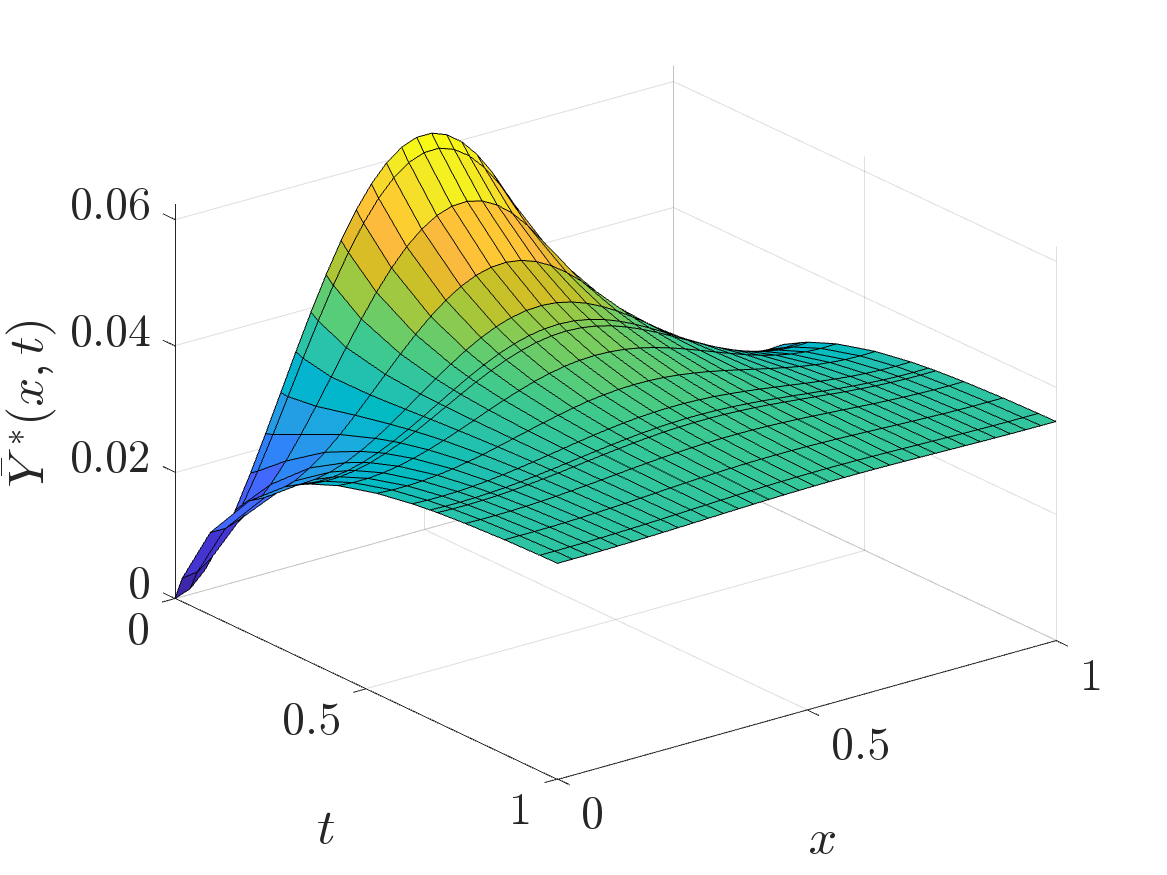}}
    \hfill
    \subfloat[Optimal controls\label{fig:burgerscontrols}]{\includegraphics[scale=0.4]{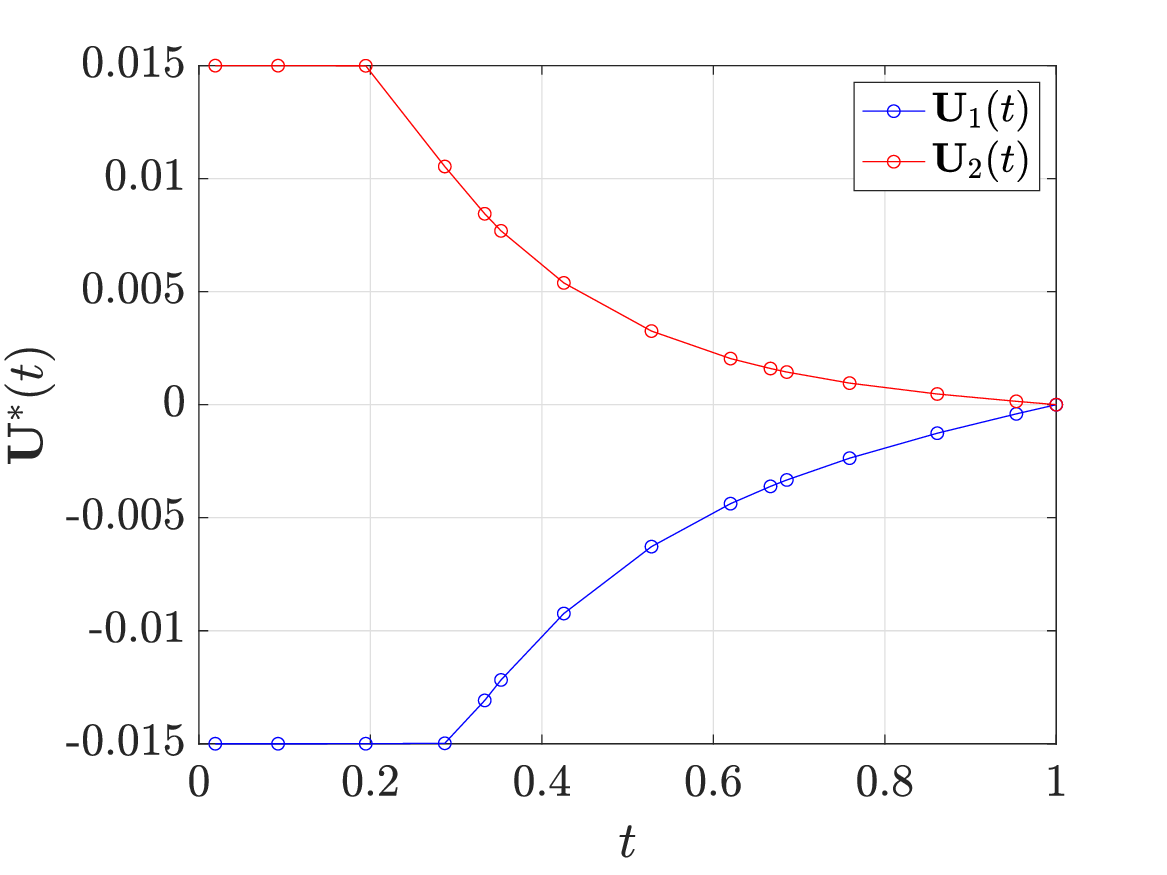}}
    \caption{The optimal state and optimal controls for the Burgers' equation tracking problem with $J=3$, $N_t = 5$, and $N_x = 34$. \label{fig:BurgersSoln}}
\end{figure}

The Burgers' equation tracking problem presented by Ref.~\cite{BuskensGriesse2006} and later solved in Ref.~\cite{Betts2020} is addressed. The problem is to minimize 
\begin{equation}
\mathcal{J} = \frac{1}{2}\int_0^1 \int_0^1 [y(x,t)-0.035]^2 \:\mathrm{d}x\mathrm{d}t + \frac{\sigma}{2} \int_0^1 [u_1^2(t)+u_2^2(t)]\:\mathrm{d}t,
\end{equation}
subject to the PDE
\begin{equation}
\frac{\partial y(x,t)}{\partial t} = \nu \frac{\partial^2 y(x,t)}{\partial x^2} - y(x,t)\frac{\partial y(x,t)}{\partial x},
\end{equation}
with the initial conditions
\begin{equation}
y(x,0) = x^2(1-x)^2,
\end{equation}
the Neumann boundary conditions
\begin{align}
\frac{\partial y(0,t)}{\partial x} &= u_1(t), \\
\frac{\partial y(1,t)}{\partial x} & = u_2(t),
\end{align}
and control bounds for $i = 1,\:2$:
\begin{equation}
u_{\text{min}} \leq u_i(t) \leq u_{\text{max}}.
\end{equation}
The problem is defined on the domain 
\begin{equation}
\Omega = \{(x,t)\:\:|\:\:0\leq x \leq 1;\:\: 0\leq t \leq 1\},
\end{equation}
and the following parameters complete the definition:
$$
\sigma = 0.01, \:\:\:\:\:\nu = 0.1,\:\:\:\:\:u_{\text{max}} = 0.015,\:\:\:\:\:u_{\text{min}} = -0.015.
$$
The NLP constraint Jacobian for the Burgers' equation tracking problem is provided in Fig.~\ref{fig:constrJacobBurgers} for a mesh with $N_t = 3$ collocation points per interval in $J = 2$ time intervals and $N_x = 5$ elements. It is noted that this mesh is too coarse for a valuable solution, but it is provided for insight and comparison with the structure presented in Fig.~\ref{fig:JacobianStructure}.
The problem is solved using the presented framework with P-1 Lagrange elements on four meshes. The computed solutions for the optimal state and the two control inputs are provided in Fig.~\ref{fig:BurgersSoln}, and the solution details are provided in Table~\ref{tab:Burgers}.



\begin{table}[hbt!]
\caption{\label{tab:Burgers}Solution details for the Burgers' equation tracking problem}
\centering
\begin{tabular}{lcccccc}
\hline
Mesh & $N_t$ & $J$ & $\mathcal{N}_t+1$ & $N_x$ & $\mathcal{J}^*$ & CPU Time (s) \\\hline
1 & 5 & 3 & 16 & 34 & $2.8709506 \times 10^{-5}$ & 2.1124  \\
2 & 3 & 10 & 31 & 34 & $2.8709897 \times 10^{-5}$  & 2.1665 \\
3 & 5 & 3 & 16 & 68 & $2.8905775 \times 10^{-5}$ & 2.2266 \\
4 & 2 & 22 & 45 & 68 & $2.8903518 \times 10^{-5}$ & 2.5911 \\
Ref.~\cite{Betts2020} & - & - & 30 & 34 & $2.89305571 \times 10^{-5}$ & 0.5020  \\
Ref.~\cite{Betts2020} & - & - & 45 & 68 & $2.89680755 \times 10^{-5}$ & 2.660  \\
\hline
\end{tabular}
\end{table}

The solutions found on each mesh were computed with an NLP tolerance of $10^{-10}$. The number of temporal collocation points in each interval was fixed. Note that the solutions on Mesh 1 and Mesh 3 demonstrate a significant reduction in the required number of temporal points to compute an accurate value of the optimal objective (as computed by another method). The solutions closely match the objective values computed in Ref.~\cite{Betts2020}. For comparison, it is noted that the solutions obtained in Ref.~\cite{Betts2020} are obtained via a temporal mesh-refinement procedure with a trapezoidal method and Hermite-Simpson separated scheme. The spatial discretization is a second-order finite difference scheme. In this work, CPU time is computed by an average over 10 trials and is defined by the total time required to execute the solution procedure (including \textbf{adigator} derivative file generation and NLP solution). It is noted that the inclusion of CPU time is meant to be used as a baseline for comparison with existing and future approaches.

As no analytical solution exists to the Burgers' equation tracking problem, a self-convergence analysis is performed for the temporal and spatial discretizations. In a self-convergence study, a solution computed on a fine mesh is treated as a reference solution, where additional solutions on coarse meshes are compared and the error can be analyzed. The analysis in Ref.~\cite{HagerHou2019} predicts \textit{at least} a convergence rate of $\mathcal{O}(h^{N_t -1})$ at the collocation points in the ``sup-norm'' (or $L_\infty-$norm) for the temporal discretization for ODE-constrained optimal control problems. However, for a self-convergence analysis (i.e.~in the absence of an analytical solution), it is infeasible to ``measure'' the error at the discrete collocation points, as the collocation points are unevenly spaced, and constructing sequences of coarse meshes that contain the points of the fine mesh is an impractical task. As a result, we depend upon the state Lagrange polynomial interpolant to interpolate the state solution to a common set of points from which we may directly compare the solution data. 

In this paper, we consider the convergence of the temporal method at a single point on the spatial grid, and likewise, we consider the convergence of the spatial method at a single point on the temporal grid. To compare the convergence rates with predicted rates in the literature, we employ the $L_{\infty}-$norm for the temporal convergence analysis \cite{HagerHou2019} and the $L_2-$norm \cite{BrennerRidgway2008} for the spatial convergence analysis.  For our purposes, we define the $L_{\infty}-$error in the discrete sense by 
\begin{equation}
    {L_{\infty} \: \mathrm{Error}} = \max_{i = 1}^N \left|\bar{Y}_i^{fine}(x_c) -\bar{Y}_i^{coarse}(x_c)\right|, 
\end{equation}
where $\bar{Y}(x_c)$ indicates the temporal state approximation at a \textit{single} spatial point, $x_c$, and $N$ is the total number of points at which the state approximations are compared. In this work, the total number of comparison points is taken to be $2\mathcal{N}_t$, with $2N_t$ comparison points in each interval, where $N_t$ is held fixed across each mesh interval. Additionally, we define a relative $L_2-$error in the discrete sense by 
\begin{equation}
    {L_{2} \: \mathrm{Error}} = \sqrt{ \dfrac{1}{N_x} \sum_{i = 1}^{N_x} \left(\dfrac{y^{fine}_{h} (x_i^{coarse},t_c) - y^{coarse}_{h}(x_i^{coarse},t_c) }{y^{fine}_{h} (x_i^{coarse},t_c) } \right)^2}. \label{eq:L2error}
\end{equation}
Here, $t_c$ represents the point in time at which the spatial state approximation is evaluated, and $N_x$ is, in this case, the number of points in the coarse solution. As the standard Galerkin finite element method employs equidistant gridpoints, it is possible to directly evaluate the $L_2-$error at the collocation points (namely, all points in the coarse solution) and an interpolant is not required (although it is included in Eq.~(\ref{eq:L2error}) for completeness). It is noted that both $t_c$ and $x_c$ exist as support points on the discretized domain. 

\begin{figure}[h!]
    \subfloat[temporal error vs. interval width \label{fig:Burgerstempconverge}]{\includegraphics[scale=0.45]{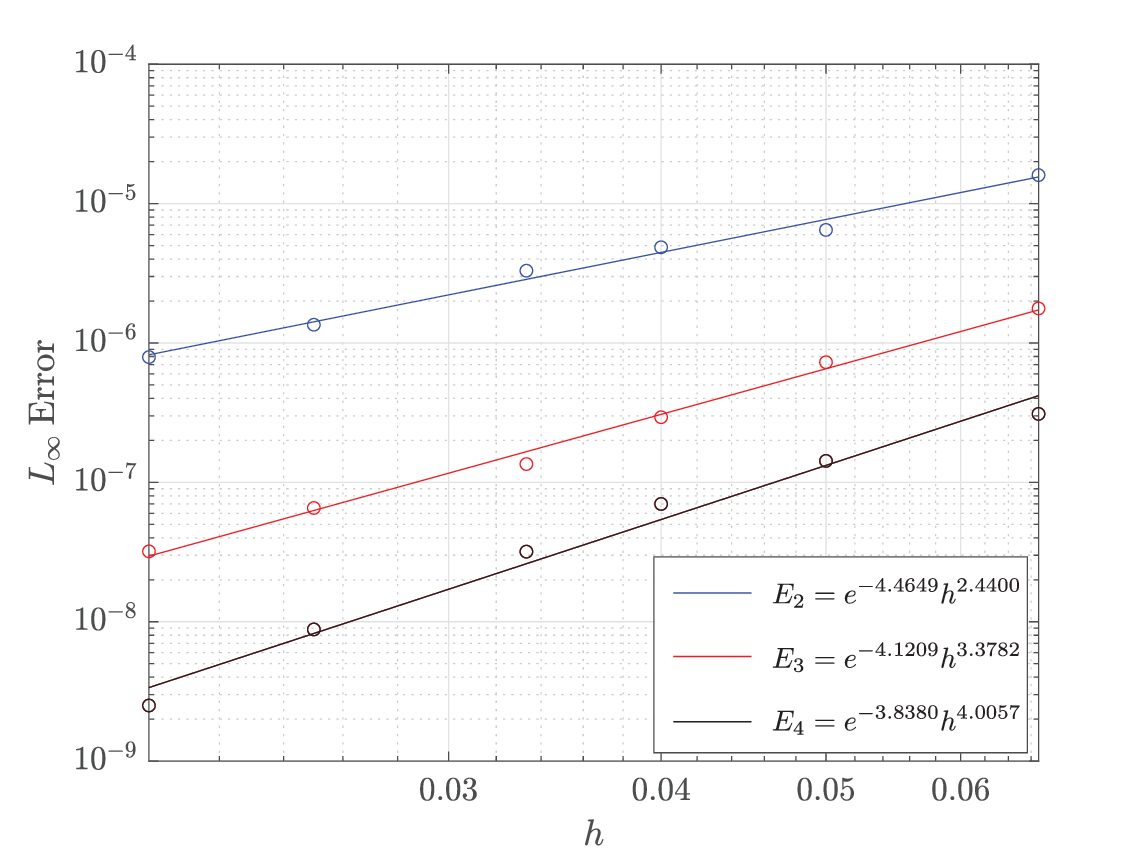}}
    \hfill
    \subfloat[spatial error vs. element width\label{fig:burgersspatconverge}]{\includegraphics[scale=0.45]{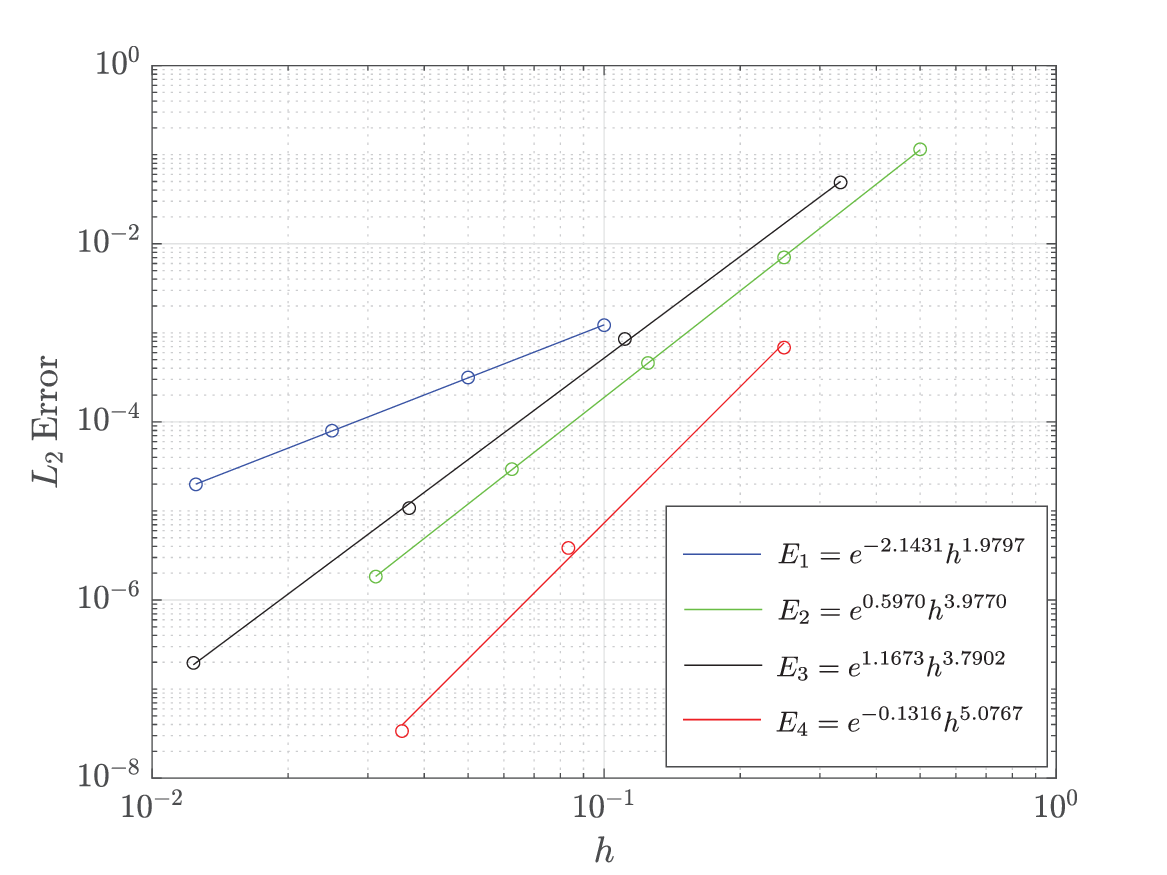}}
    \caption{Temporal and spatial point-convergence of the state solutions at $x_c = 0.2388$ (a) and $t_c = 1$ (b), respectively. \label{fig:convBurgers}}
\end{figure}

For the Burgers' equation problem, the temporal point used for comparison was selected as $t_c = 1$ to allow for sufficient time relaxation and the spatial point used for comparison was selected as $x_c = 0.2388$. The convergence rates in the $L_\infty-$ and $L_2-$senses are demonstrated in log-log plots of the error vs.~the mesh size, $h$, in Fig.~\ref{fig:convBurgers}. For the temporal convergence study, $h$ represents the interval width, and for the spatial convergence study, $h$ represents the element width. Additionally, in Fig.~\ref{fig:convBurgers}, the subscript $i$ in $E_i$ indicates the degree of the state approximation in each interval/element. 

In Fig.~\ref{fig:Burgerstempconverge}, the convergence rate roughly appears at $\mathcal{O}(h^{N_t})$ (or slightly faster), which is faster than the conservative bound presented in Ref.~\cite{HagerHou2019}, but slower than the empirically computed result of $\mathcal{O}(h^{N_t+2})$ in Ref.~\cite{HagerHou2019}. Likely, the deviation in the apparent convergence rate is due to the performance of a self-convergence study and the evaluation of the error at points that are \textit{not} the collocation points. In Fig.~\ref{fig:burgersspatconverge}, the convergence rate is roughly $\mathcal{O}(h^{p+1})$, where $p$ indicates the degree of the state approximation in each element, which is expected for piecewise continuous Galerkin finite elements. For $p = 2$ elements, a superconvergence result is demonstrated where the error decays at a rate of approximately $\mathcal{O}(h^4)$. This result is likely due to sufficient solution regularity and the evaluation of the error at the nodal points \cite{wahlbin2006,Douglas1974}. In both cases, the error converges exponentially as a function of the mesh size.

\subsection{Heat Equation Problem}
An optimal control problem described by the heat equation is presented. It can be viewed as a simplified model for the heating of a probe in a kiln. The problem is solved in Refs.~\cite{Betts2020,Heinkenschloss1996}. The goal is to minimize the deviation from a desired temperature profile as defined by the objective
\begin{equation}
    \mathcal{J} = \frac{1}{2}\int_0^{t_f} \left\{[y(1,t)-y_d(t)]^2+\gamma u^2(t)\right\}\:\mathrm{d}t,
\end{equation}
by choosing the control function between allowable limits
\begin{equation}
    u_{\min} \leq u(t) \leq u_{\max}, \label{eq:controllim}
\end{equation}
that satisfies the nonlinear heat equation
\begin{equation}
        q(x,t) = (a_1+a_2y)\frac{\partial y}{\partial t}-a_3\frac{\partial^2 y}{\partial x^2}-a_4\left(\frac{\partial y}{\partial x}\right)^2-a_4y\frac{\partial^2 y}{\partial x^2},
\end{equation}
and the initial and boundary conditions 
\begin{align}
        (a_3+a_4y)\left. \frac{\partial y}{\partial x}\right|_{x=0} &= g[y(0,t)-u(t)], \\
    (a_3+a_4y)\left. \frac{\partial y}{\partial x}\right|_{x=1} &= 0, \\
    y(x,0) &= y_{I}(x),
\end{align}
where the following definitions: 
    \begin{align*}
        y_d(t) =&\:2-\mathrm{e}^{\rho t}, \\
        y_I(x) =&\:2 + \mathrm{cos}(\pi x), \\
        q(x,t) =&\:[\rho(a_1+2a_2)+\pi^2(a_3+2a_4)]\mathrm{e}^{\rho t}\mathrm{cos}(\pi x), \\ 
        &-a_4\pi^2\mathrm{e}^{2\rho t}+(2a_4\pi^2+\rho a_2)\mathrm{e}^{2\rho t}\mathrm{cos}^2(\pi x),
    \end{align*}
    and parameters in Eq.~\eqref{eq:neccparameters} complete the problem statement:
\begin{equation}
    \begin{array}{lclclcl}
        a_1 & = & 4 &, & \rho & = & -1, \\
        a_2 & = & 1 & , & t_f & = & 0.5, \\ 
        a_3 & = & 4 & , & \gamma & = & 10^{-3}, \\
        a_4 & = & -1 & , & g & = & 1, \\ 
        u_{\min} & = & -\infty & , & u_{\max} & = & 0.1.
    \end{array}
    \label{eq:neccparameters}
\end{equation}

\begin{figure}[h]
    \centering
    \includegraphics[scale = 0.325]{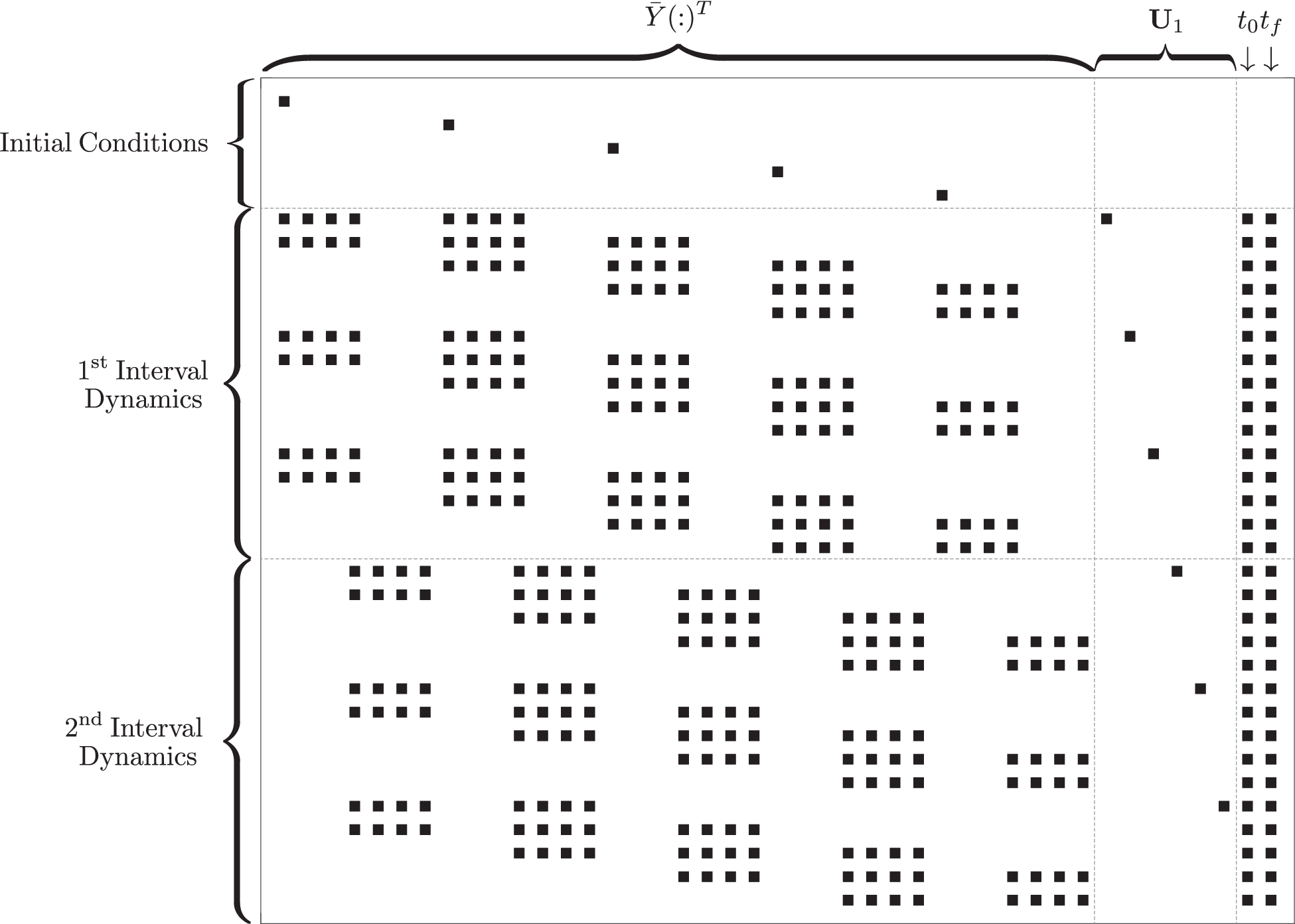}
    \caption{NLP constraint Jacobian for the heat equation problem with P-1 Lagrange elements and $N_t^{(1)} = N_t^{(2)} = 3$, $J = 2$, $N_x = 5$.}
    \label{fig:heatJacob}
\end{figure}

\begin{figure}[h]
    \subfloat[Optimal state \label{fig:heatstate}]{\includegraphics[scale=0.4]{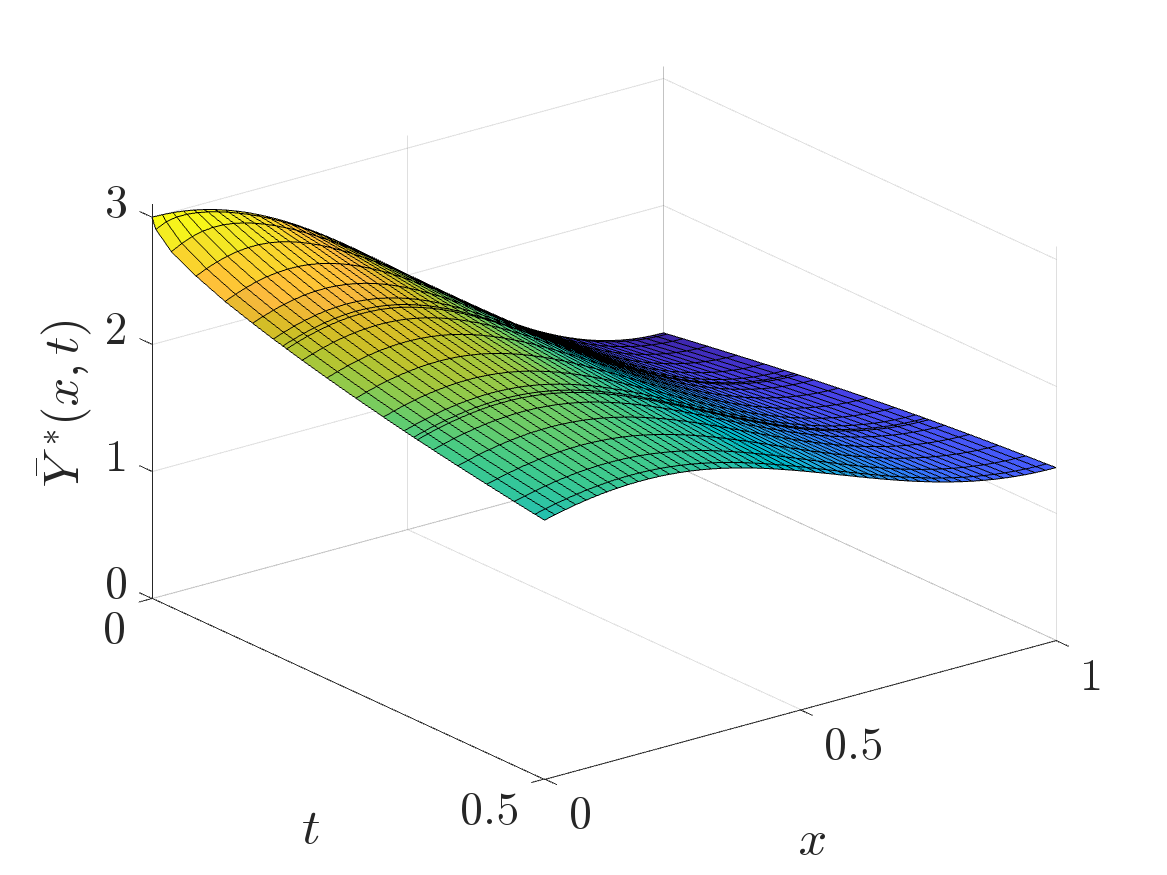}}
    \hfill
    \subfloat[Optimal control \label{fig:heatcontrol}]{\includegraphics[scale=0.4]{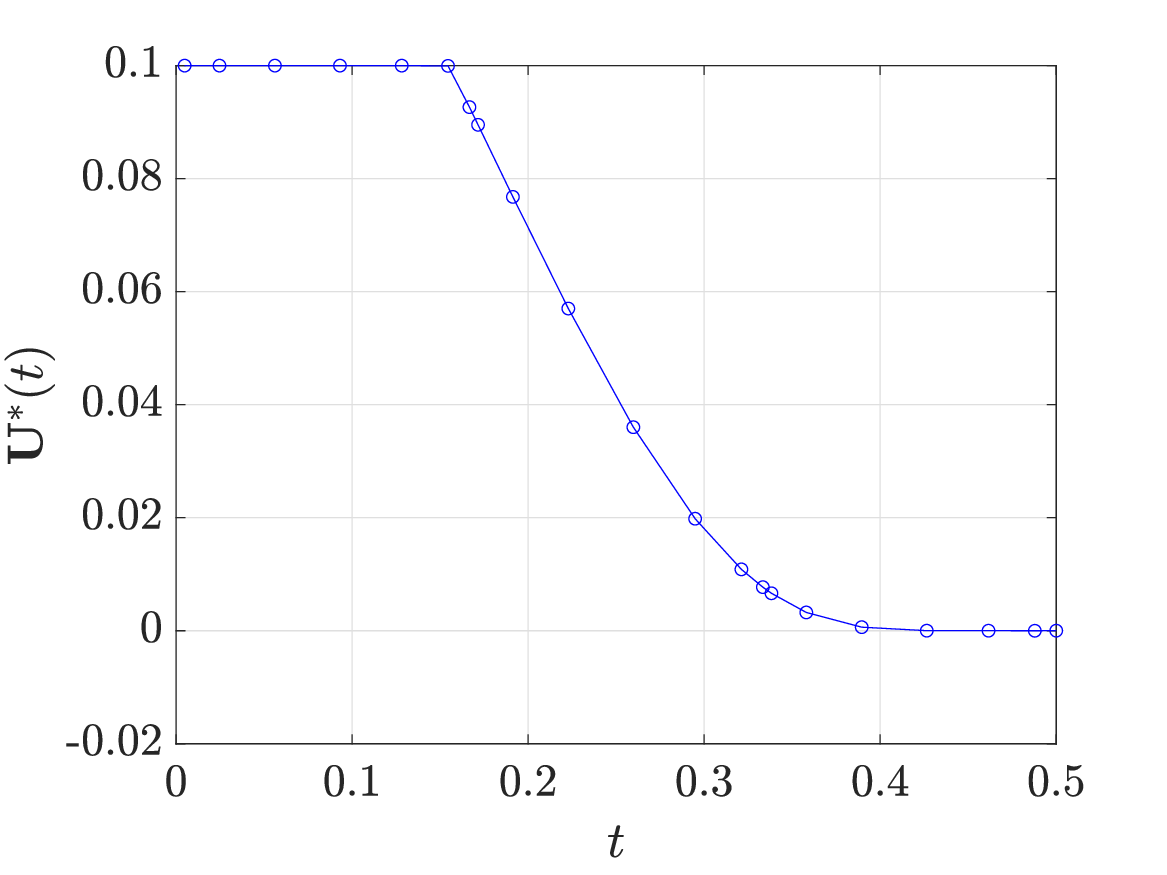}}
    \caption{The optimal state and optimal control for the heat equation problem with $J=3$, $N_t = 7$, and $N_x = 50$. \label{fig:HeatSoln}}
\end{figure}

   
In this example, the state $y(x,t)$ is the temperature, the functions ($a_1 + a_2 y$) and ($a_3 + a_4y$) represent the temperature-dependent specific heat capacity and thermal conductivity, respectively, $q(x,t)$ is a source term, $y_d(t)$ is a given desired temperature profile, $y_I(x)$ is the initial temperature profile, $\gamma$ is a positive regularization parameter, and $g$ and $\rho$ are given scalars. An example NLP constraint Jacobian for the heat equation problem is provided in Fig.~\ref{fig:heatJacob} for the same mesh used in Fig.~\ref{fig:constrJacobBurgers}. It is evident that due to the existence of only one near boundary control, the portion of the constraint Jacobian corresponding to $\nabla_{\mathbf{U}_2}\mathbf{F}$ in Fig.~\ref{fig:JacobianStructure} is left absent. 
The problem is solved using the presented framework with P-1 Lagrange elements on 4 meshes. The computed solutions for the optimal state and control input are provided in Fig.~\ref{fig:HeatSoln}, and the solution details are provided in Table \ref{tab:HeatFRFEComp}.

\begin{table}[h]
\caption{\label{tab:HeatFRFEComp}Solution details for the heat equation problem}
\centering
\begin{tabular}{lcccccc}
\hline
Mesh & $N_t$ & $J$ & $\mathcal{N}_t+1$ & $N_x$ & $\mathcal{J}^*$ & CPU Time (s) \\\hline
1 & 7 & 3 & 22 & 20 & $3.6232288 \times 10^{-5}$  & 2.3989 \\
2 & 3 & 33 & 100 & 50 & $3.8283815 \times 10^{-5}$ & 3.3314 \\
3 & 7 & 3 & 22 & 50 &  $3.8283491 \times 10^{-5}$ & 2.4886 \\
4 & 4 & 10 & 41 & 50 & $ 3.8283552  \times 10^{-5}$ & 2.6656 \\
Ref.~\cite{Betts2020} & - & - & 99 & 50 & $3.8786845 \times 10^{-5}$ & 8.86  \\
Ref.~\cite{Heinkenschloss1996} & - & - & 100 & 20 & $4.03 \times 10^{-5}$ & -  \\
\hline
\end{tabular}
\end{table}
The solutions found on each mesh were computed with an NLP tolerance of $10^{-10}$. The number of temporal collocation points in each interval was fixed. Similar to the results of the viscous Burgers' problem, it is demonstrated that comparable values of the optimal objective can be computed with fewer temporal points. For comparison, the discretization in Ref.~\cite{Heinkenschloss1996} is a piecewise constant temporal discretization with P-1 linear finite elements, and the discretization in Ref.~\cite{Betts2020} is a sixth-order Lobatto IIIA temporal discretization with second-order finite differencing in the spatial dimension.

\begin{figure}[h!]
    \subfloat[temporal error vs. interval width \label{fig:Heattempconverge}]{\includegraphics[scale=0.45]{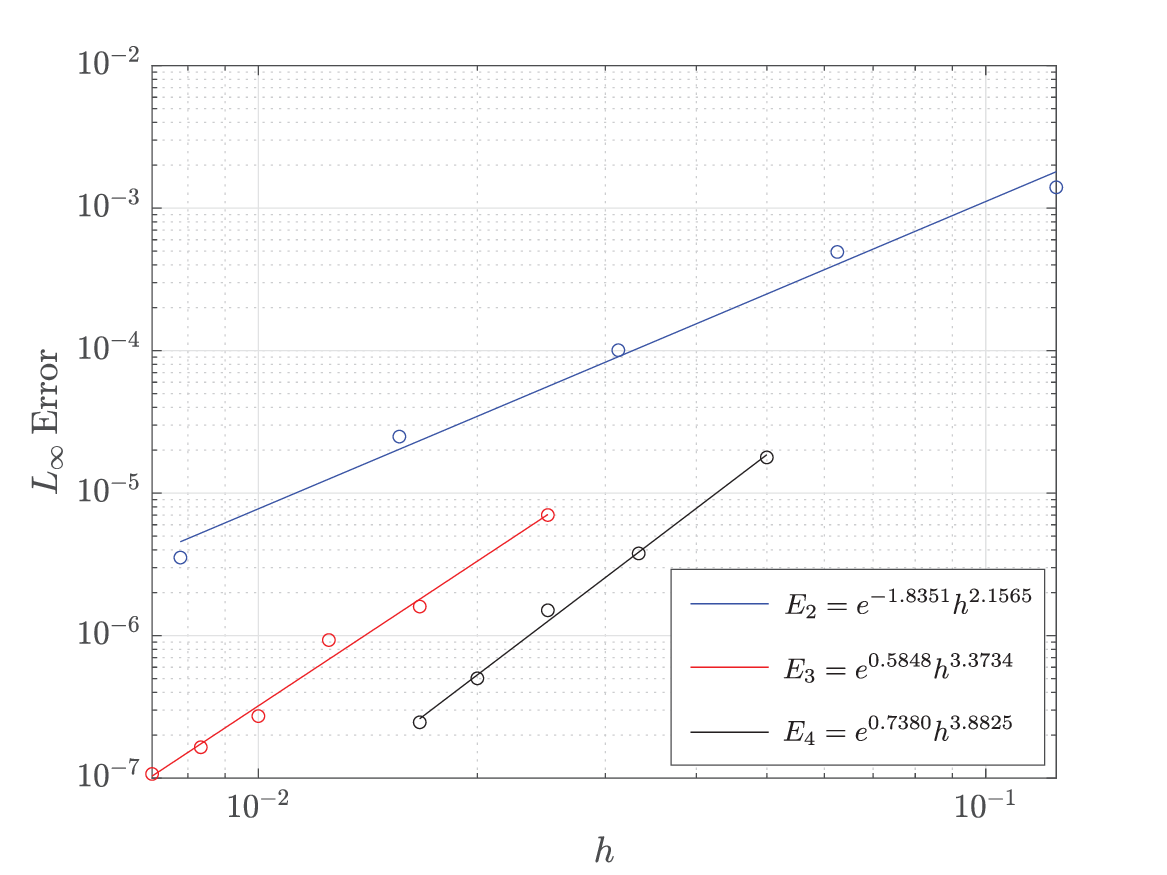}}
    \hfill
    \subfloat[spatial error vs. element width\label{fig:heatspatconverge}]{\includegraphics[scale=0.45]{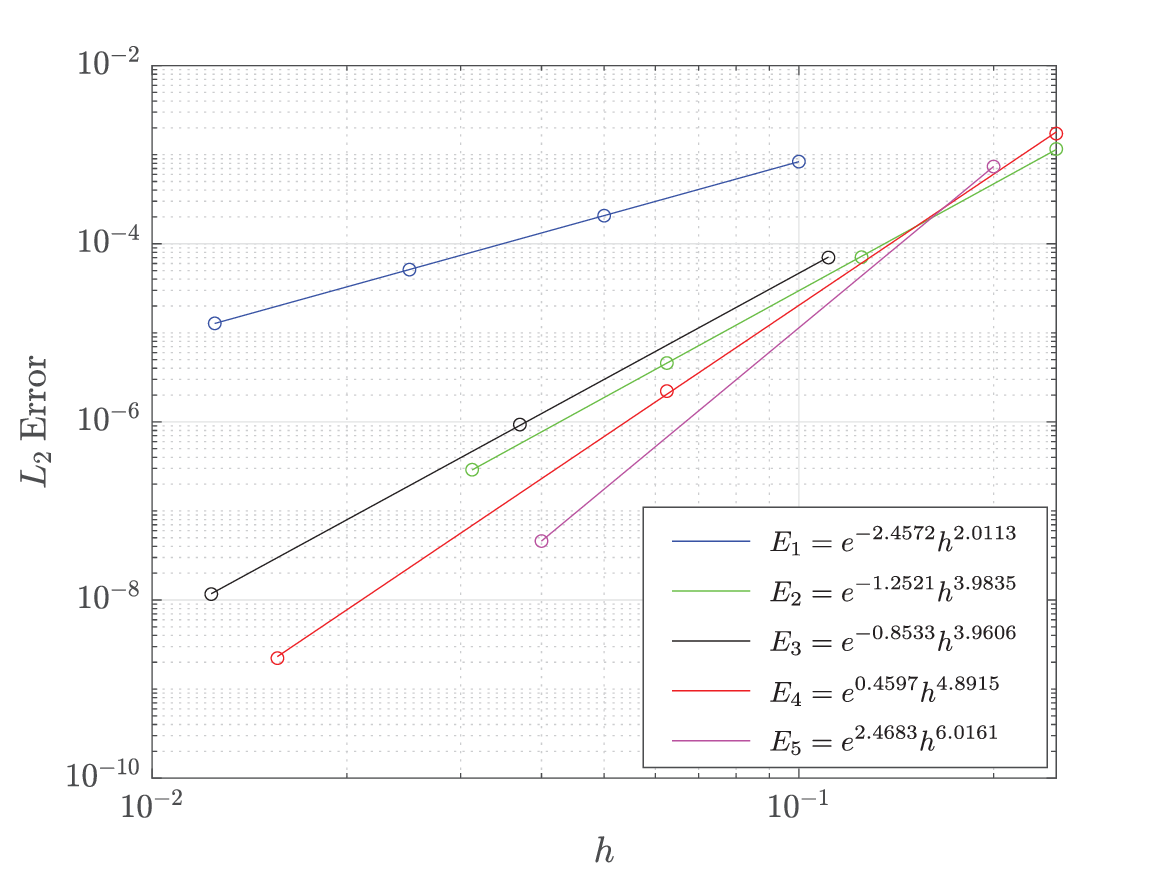}}
    \caption{ Temporal and spatial point-convergence of the state solutions at $x_c = 1/3$ (a) and $t_c = 0.5$ (b), respectively.} \label{fig:convHeat}
\end{figure}

Similarly to the Burgers' equation problem, no analytical solution exists to the heat equation problem. As a result, we perform a self-convergence analysis to demonstrate the convergence of the method on an additional example. Similarly to the Burgers' equation problem, we select $t_c = 0.5$ (or the endpoint of the time domain) to allow for sufficient relaxation. For the temporal convergence analysis, we select $x_c = 1/3$. The temporal and spatial convergence rates in the $L_\infty-$ and $L_2-$senses are demonstrated in log-log plots of the error vs.~the mesh size in Fig.~\ref{fig:convHeat}.
As is demonstrated in Fig.~\ref{fig:convHeat}, the error decreases exponentially as a function of the mesh size in both dimensions. The temporal error decays at roughly $\mathcal{O}(h^{N_t})$, and the spatial error decays at roughly $\mathcal{O}(h^{p+1})$. Similarly to the Burgers' equation case, a superconvergence result is demonstrated for P-2 elements due to the regularity of the solution at the final time.

\subsection{Heat Equation Problem with Control Constraints}

\begin{figure}[h]
    \subfloat[Optimal state \label{fig:heatstateconst}]{\includegraphics[scale=0.4]{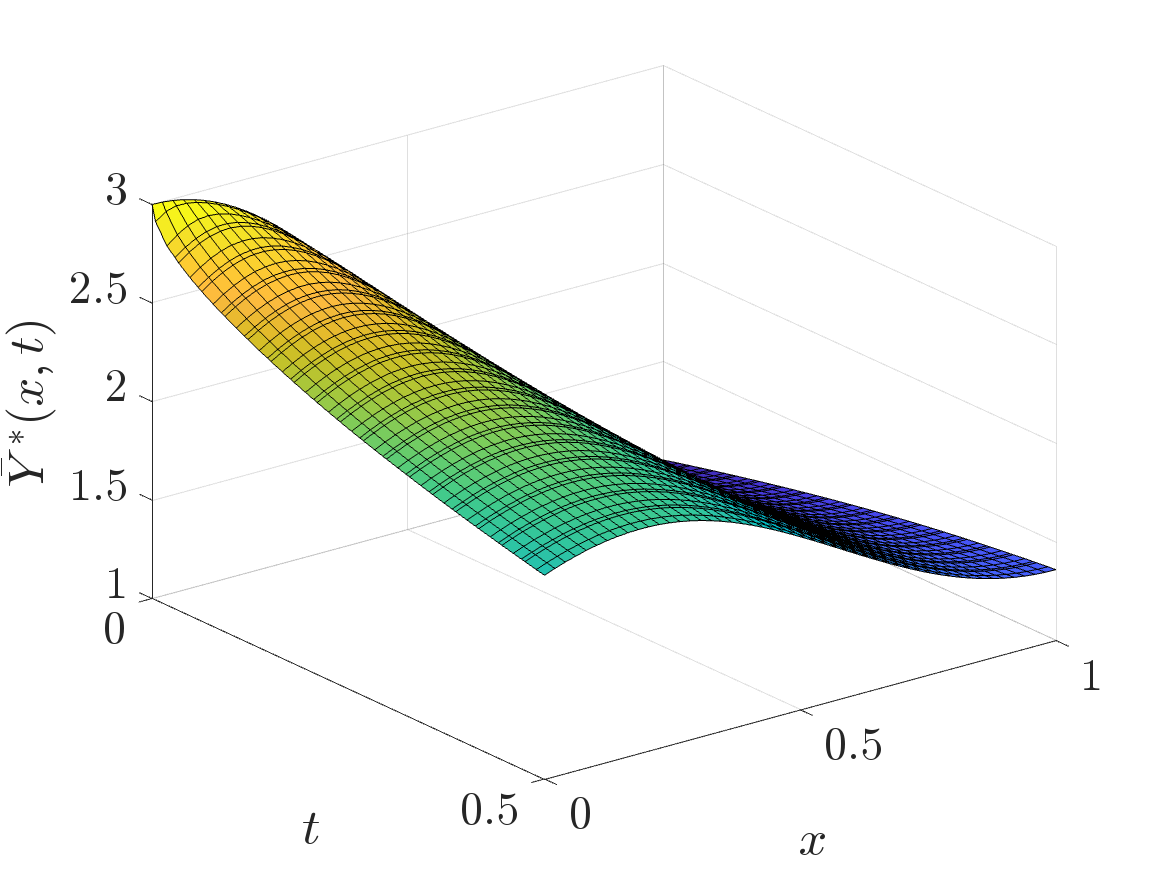}}
    \hfill
    \subfloat[Optimal control \label{fig:heatcontrolconst}]{\includegraphics[scale=0.4]{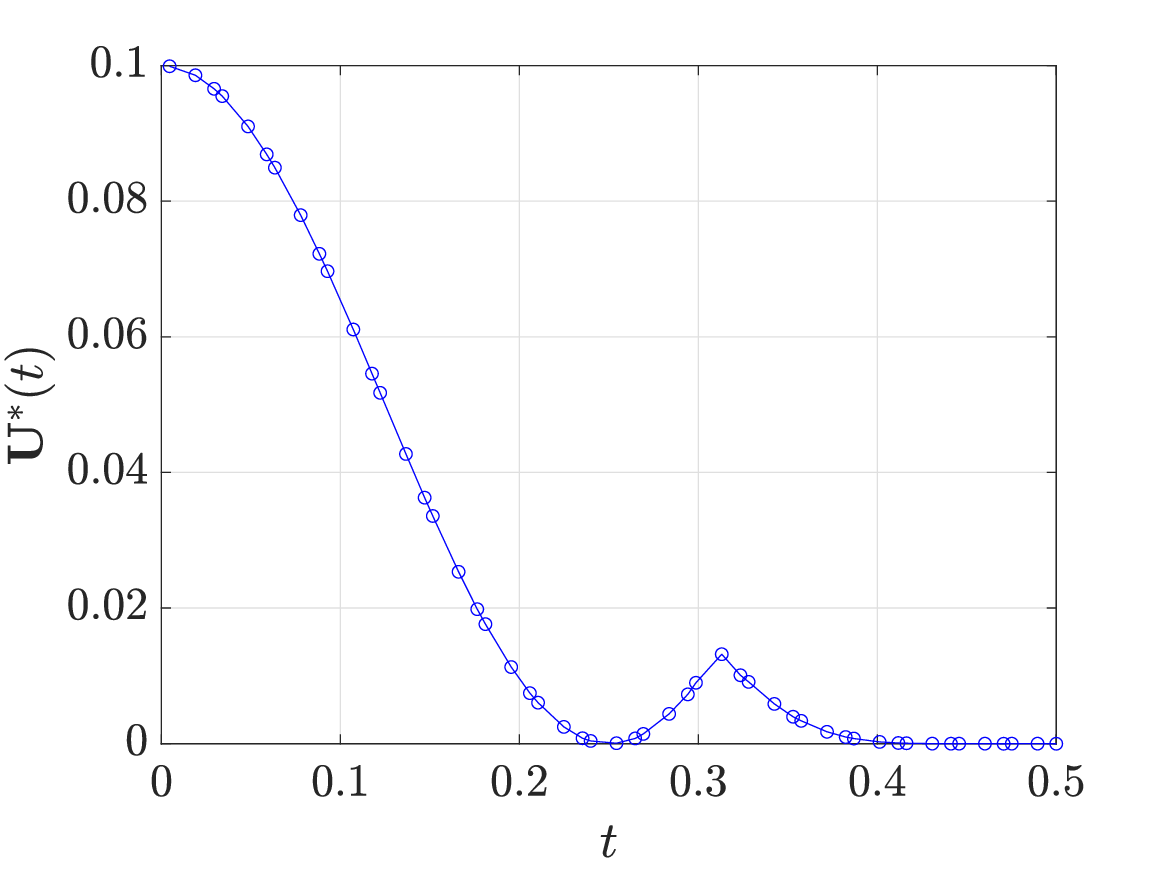}}
    \caption{ The optimal state and optimal control for the heat equation problem with $J=17$, $N_t = 3$, and $N_x = 50$.} \label{fig:newconst}
\end{figure}

Another benefit of the method presented in this paper is that we may easily apply time-dependent inequality constraints in the problem formulation. If we apply an additional constraint to Eq.~(\ref{eq:controllim}) such that 
\begin{equation}
    u(t) \leq u_{\max} \dfrac{1+\mathrm{cos}(4 \pi t)}{2},
\end{equation}
the constraint can be enforced implicitly through the upper limits of the control variable specified in the NLP. The problem with the additional control constraint is solved using the presented framework with P-1 Lagrange elements on 4 meshes. The computed solutions for the optimal state and control on Mesh 1 from Table \ref{tab:HeatNewConst} are provided in Fig.~\ref{fig:newconst}. Notably, the control structure lies along the inequality constraint until shortly after $t = 0.3$, where, after, the control relaxes to zero. Due to the additional constraint, the value of the optimal objective slightly increases, which implies additional cost to meet the problem objective. The use of a direct framework allows for the user to implement inequality constraints with relative ease, which further motivates the method presented in this work. 

\begin{table}[h]
\caption{\label{tab:HeatNewConst} Solution details for the heat equation problem with control constraints}
\centering
\begin{tabular}{lcccccc}
\hline 
Mesh & $N_t$ & $J$ & $\mathcal{N}_t+1$ & $N_x$ & $\mathcal{J}^*$ & CPU Time (s) \\\hline
1 & 3 & 17 & 52 & 50 & $3.8669419 \times 10^{-5}$ & 2.7897 \\
2 & 5 & 10 & 51 & 50 & $3.8669506 \times 10^{-5}$ & 2.7233 \\
3 & 3 & 17 & 52 & 100 & $3.8954568 \times 10^{-5}$ & 3.0697 \\
4 & 5 & 10 & 51 & 100 & $ 3.8954649  \times 10^{-5}$ & 3.3684 \\
\hline
\end{tabular}
\end{table}

\subsection{Discussion of Results}
 As is common for optimal control problems governed by PDEs, there is no analytical solution to the presented problems. So, in place of an analytical solution, a numerical comparison to Ref.~\cite{Betts2020} was performed. To compare the optimal objective on a comparable mesh, solutions using the same number of points as Ref.~\cite{Betts2020} were desired. However, to produce a mesh with $\mathcal{N}_t+1 = 30$ points requires a dynamic mesh (one in which the number of collocation points may vary per interval). This is a motivation for future implementation. In this work, solutions are computed on a static mesh with a fixed number of collocation points per interval. As a result, a mesh with the exact same number of points as Ref.~\cite{Betts2020} is, at the moment, unattainable. To get as close as possible, Mesh 2 in Table \ref{tab:Burgers} was constructed. It is demonstrated that the value of the computed objective is slightly smaller than the result computed by Ref.~\cite{Betts2020} in a larger CPU time. The CPU time included in Table \ref{tab:Burgers} is comparably larger than the solution obtained in Ref.~\cite{Betts2020} (for meshes of smaller size) as it contains the time required to compute derivatives via automatic differentiation \textit{and} the time required to produce a solution to the NLP. Generally, the overhead associated with automatic differentiation is the limiting factor in the total CPU cost. One could mitigate the CPU cost through analytical derivatives, or derivative approximation techniques such as sparse central-differencing, bicomplex-step, or hyper-dual methods \cite{AgamawiRao2020}, but these methods are beyond the scope of this manuscript. The solution on Mesh 1 was provided to demonstrate that the variation in the computed objective is small subject to drastic decreases in the number of temporal points. The solutions on Mesh 3 and Mesh 4 further demonstrate this property. 

For the heat equation problem, a similar effect is demonstrated. The solution on Mesh 2 contains the same number of points as the solution in Ref.~\cite{Heinkenschloss1996}, wherein a piecewise constant temporal integration scheme was used (and is likely the reason for the relatively large value of the computed optimal objective). The solution in Ref.~\cite{Betts2020} demonstrates a closer value of the computed optimal objective to our solutions. It is evident that on Mesh 1, too few spatial elements are used to provide a computed objective that resembles the solutions on the other meshes. The solutions on Mesh 3 and Mesh 4 demonstrate that the optimal objective is largely unaffected by the reduction in the number of temporal points, which highlights that fewer numbers of temporal points may be used to compute accurate solutions. For this problem, where the solution structure is relatively smooth, a Lagrange polynomial approximation of the state with fLGR support points provides a good underlying approximation of the solution structure, which provides a justification for this phenomenon. 

This leads us to conclude that this method is well suited for problems in which the solution structure in time is relatively smooth. When this is the case, a polynomial approximation of the state is suitable for the solution structure. In this work, the number of collocation points per interval and number of intervals were largely selected based on observation. The number of intervals and number of collocation points per interval that provide the best solutions are those that possess mesh points as close as possible to drastic changes in the control structure. This interchange is similar to the findings of Ref.~\cite{DarbyHager2011a}, where it is emphasized that regions with rapid changes in the solution structure should be presented with fewer collocation points per interval and a larger number of intervals. To improve upon this work, it is desired that a dynamic mesh be developed in which the locations of the mesh intervals and collocation points be determined algorithmically based on the solution structure.     

On the topic of CPU time, in comparison to the results in Ref.~\cite{Betts2020}, the computational speed is slower for the Burgers' equation problem for meshes of smaller size and faster for the heat equation problem. The CPU time could likely be improved through the NLP solver, automatic differentiation software, code implementation, or simply an improved initial guess, where, in this work, a straight-line guess from the initial condition to zero is used. In general, longer CPU times in this work are attributed to time spent generating the \textbf{adigator} derivative files for {\em IPOPT} and the NLP solver itself. For example, for the solution on Mesh 4 in Table~\ref{tab:Burgers}, approximately 2.1779 seconds were spent generating the derivative files for {\em IPOPT}, and approximately 0.4034 seconds were spent within {\em IPOPT} itself (where both times are computed as averages over 10 trials). Thus, to improve CPU performance, one could address the NLP derivative approximation technique to seek the most improvement.


\section{Conclusions}\label{sec:conclusions}
 A framework for the numerical solution of optimal control problems governed by parabolic PDEs has been presented. A multi-interval fLGR collocation method was paired with a finite element spatial discretization that allowed for the satisfaction of boundary/dynamics relationships at the terminal time. In the optimal control of ODEs, it has been demonstrated that the fLGR and standard LGR temporal discretizations may be used nearly interchangeably to produce comparable results. It was demonstrated that for the optimal control of systems governed by PDEs, the fLGR discretization is preferable. The use of the finite element spatial discretization further allowed for a proper NLP search space and a reduction in the problem size in comparison to the method presented in Ref.~\cite{DaviesDennis2024}, where an introduction of additional state variables and a user redefinition of constraints was required. In addition, to handle nonlinearities present in the variational form, a generalization of a Kirchoff transformation was 
 used to pose the variational form as a linear system in an additional variable. The procedure to transcribe the optimal control problem into a nonlinear programming problem was outlined, and derivatives of the NLP objective and constraint functions were summarized. The presented method produced comparable values of the optimal objective to other methods on two test problems. It was shown on each test problem that multi-interval orthogonal collocation can lead to a reduction in the required number of temporal points to compute accurate values of the objective. Finally, numerical results demonstrated an exponential decay of the error with respect to the mesh size in both the temporal and spatial dimensions.

\section*{Funding}

{\noindent}The authors gratefully acknowledge support for this research from the U.S. Air Force Research Laboratory under grant FA8651-25-1-0002, from the U.S. National Science Foundation under grant CMMI-2031213, and from the U.S. Office of Naval Research under grant N00014-22-1-2397.

\section*{Acknowledgments}

{\noindent}The authors acknowledge John T.~Betts for providing his solution to the heat equation problem found in Section \ref{sec:NumericalExamples} of this paper.  

\section*{Disclaimer}

{\noindent}The views and conclusions contained herein are those of the authors and should not be interpreted as necessarily representing the official policies or endorsements, either expressed or implied, of AFRL/RW or the U.S. Government.

\appendix  
  \renewcommand{\theequation}{A-\arabic{equation}}
  \setcounter{equation}{0}  
\section*{Appendix}
The framework introduced is not restricted to a spatial discretization with finite elements. For implementation with finite differencing, the procedure utilizes the method of lines. For further detail, one can refer to Ref.~\cite{DaviesDennis2025}. Similar to finite elements, the problem is semi-discretized in space prior to the temporal discretization. The semi-discretization procedure forms a coupled differential-algebraic system that can be discretized in time with a well-suited collocation scheme to accurately capture the effects of the boundary conditions. As opposed to the finite element method, the finite difference method is utilized to approximate the strong form of the PDE. Starting from Eq.~\eqref{eq:strongpde}, one has 
\begin{equation}
         \frac{\partial y^{(j)}}{\partial r} +  \psi^{(j)} \kappa(y^{(j)}) \frac{\partial y^{(j)}}{\partial x} = a \psi^{(j)} \frac{\partial^2 y^{(j)}}{\partial x^2}. \label{eq:strongpde2}
\end{equation}
The spatial domain, $x \in \Omega$, can be partitioned into $K$ intervals (or $N_x = K+1$ parallel ``lines") by $\mathcal{S}_k = [x_{k-1},x_{k}] \subset \Omega,$ $\{k=1,\ldots, K\}$, where 
\begin{equation}
    \bigcup_{k = 1}^{K} \mathcal{S}_k = \Omega,\quad \bigcap_{k=1}^{K} \mathcal{S}_k = \{x_1,\ldots,x_{K-1}\},
\end{equation}
and a state variable, $y^{(j)}_k$, $k=\{0,\ldots,K\}$, is defined at each line. At the interior lines defined by $\cap_{k=1}^K \mathcal{S}_k$, central difference approximations can be applied to approximate the spatial derivatives. For simplicity in the following discussion, we will assume the use of second-order central difference approximations at equidistant grid points; however, the procedure is not dependent upon this assumption. For reference, higher-order differencing schemes are employed in Ref.~\cite{Betts2020}. Equation \eqref{eq:strongpde} can be semi-discretized by 
\begin{equation}
            \frac{\mathrm{d} y_k^{(j)}}{\mathrm{d}r} + \psi^{(j)}\kappa(y_k^{(j)})\left( \frac{y_{k+1}^{(j)}-y_{k-1}^{(j)}}{2h} \right)= a\psi^{(j)} \left(\frac{y_{k+1}^{(j)}-2y_k^{(j)}+y_{k-1}^{(j)}}{h^2} \right), \label{eq:transformedPDEFRFD3}
\end{equation}
for $k=\{1,\ldots,K-1\}$ to form a coupled set of $N_x-2$ ordinary differential equations. To complete the spatial discretization, the boundary conditions must be addressed. Dirichlet conditions are straightforward as they imply fixed conditions on the states $y^{(j)}_0$ and $y^{(j)}_K$. For boundary conditions that are functions of the state derivative (Neumann, Robin, etc.), forward and backward difference approximations can be utilized. For example, for the conditions provided in Eq.~\eqref{eq:bc1} and Eq.~\eqref{eq:bc2}, one has 
\begin{align}
     &\left.\dfrac{\partial y^{(j)}}{\partial x}\right|_{x = x_0} \approx \frac{-y^{(j)}_{2}+4y^{(j)}_1-3y^{(j)}_0}{2h} + \mathcal{O}(h^2)\:\:\:= u^{(j)}_1(t), \\
    &\left.\dfrac{\partial y^{(j)}}{\partial x}\right|_{x = x_f} \approx \frac{y^{(j)}_{K-2}-4y^{(j)}_{K-1}+3y^{(j)}_K}{2h}+ \mathcal{O}(h^2) = u^{(j)}_2(t),
\end{align}
to provide constraints on the boundaries. Thus, a differential-algebraic system of $N_x$ coupled equations is formed. The temporal discretization is applied in a similar manner to the application in Section \ref{sec:fLGR}. The state, $y_k$, $k=\{0,\ldots,K\}$, can be approximated by a basis of Lagrange polynomials by 
\begin{equation}
y_k^{(j)} \approx Y_k^{(j)}  = \sum_{n = 0}^{N_t^{(j)}} Y_{nk}^{(j)}L_n^{(j)}(r),
\end{equation}
where the Lagrange polynomials are defined identically in Eq.~\eqref{eq:Lpoly}. Differentiating the state approximation on the interior lines and applying the expression at the fLGR points, $r_i$, $i = \{1,\ldots,N_t^{(j)}\}$ provides
\begin{equation}
    \dfrac{\mathrm{d} y_k^{(j)}(r_i)}{\mathrm{d} r} \approx  \dfrac{\mathrm{d} Y_k^{(j)}(r_i)}{\mathrm{d} r} = \sum_{n=0}^{N_t^{(j)}}  \dfrac{\mathrm{d} L^{(j)}_n(r_i)}{\mathrm{d} r} Y_{nk}^{(j)}, 
\end{equation}
or, equivalently,
\begin{equation}
    \dfrac{\mathrm{d} y_k^{(j)}(r_i)}{\mathrm{d} r} \approx  \dfrac{\mathrm{d} Y_k^{(j)}(r_i)}{\mathrm{d} r} = \sum_{n=0}^{N_t^{(j)}} D_{in}^{(j)}Y_{nk}^{(j)}.
\end{equation}
A state matrix can be constructed by 
\begin{equation}
    \bar{Y} = \begin{bmatrix}
        Y_{nk}^{(1)} & Y_{ik}^{(2)} & \ldots & Y_{ik}^{(J)}
    \end{bmatrix}^T,
\end{equation}
where $i = \{1, \ldots,N_t^{(j)}\},$ $n = \{0, \ldots,N_t^{(1)}\},$ and $k = \{0, \ldots,K\},$ that allows for the construction of differentiation matrices. The temporal differentiation matrix is identical to that of Eq.~\eqref{eq:diffmat}; whereas in contrast, the spatial differentiation matrices $D_x \in \mathbb{R}^{(N_x-2) \times N_x}$ and $D_{xx} \in \mathbb{R}^{(N_x-2) \times N_x}$ are the typical bi-diagonal and tri-diagonal matrices associated with second-order differencing for first and second derivatives, respectively. This leads to the construction of a fully-discrete form of Eq.~\eqref{eq:transformedPDEFRFD3} that is constructed by 
\begin{equation}
    D_t\bar{Y}_{1:K-1} + \boldsymbol{\alpha}^T \odot \bar{\kappa}\odot (D_x(\bar{Y}_{1:\mathcal{N}_t,:})^T)^T = a \boldsymbol{\alpha}^T \odot (D_{xx}(\bar{Y}_{1:\mathcal{N}_t,:})^T)^T,
\end{equation}
with discrete boundary conditions given by 
\begin{align}
    \bar{Y}_{1:\mathcal{N}_t,0} &= \frac{4\bar{Y}_{1:\mathcal{N}_t,1}-\bar{Y}_{1:\mathcal{N}_t,2}-2h\mathbf{U}^T_1}{3}, \\
    \bar{Y}_{1:\mathcal{N}_t,K} &= \frac{4\bar{Y}_{1:\mathcal{N}_t,K-1}-\bar{Y}_{1:\mathcal{N}_t,K-2}+2h\mathbf{U}^T_2}{3},
\end{align}
and discrete initial conditions 
\begin{equation}
    \bar{Y}_{0k} = q(x_k,t_0), \quad k =\{0,\ldots,K\}.
\end{equation}
The matrix, $\bar{\kappa} \in \mathbb{R}^{\mathcal{N}_t \times (N_x-2)}$, is constructed by the evaluation of the function at the interior lines and temporal fLGR points (i.e.~$\bar{\kappa} = \kappa(\bar{Y}_{(1:\mathcal{N}_t),(1:K-1))})$). The only remaining procedure is to discretize the objective, which is done by the trapezoidal rule and an fLGR quadrature in space and time, respectively:
\begin{multline}
   \mathcal{J} \approx  \sum_{j=1}^{J}\sum_{i=1}^{N_t^{(j)}} \psi^{(j)} w_i^{(j)} \left[\frac{h}{2} \left(\mathcal{L}(x_0,t(r_i^{(j)}),Y^{(j)}_{i0})  \right.\right.
   + \left.\mathcal{L}(x_K,t(r_i^{(j)}),Y^{(j)}_{iK}) \right) \left. \right. \\ \left. +h\sum_{k=1}^{K-1} \mathcal{L}(x_k,t(r_i^{(j)}),Y^{(j)}_{ik})  
   + \mathcal{P}(t(r_i^{(j)}),U_1(r_i^{(j)}),U_2(r_i^{(j)}),Y_{i0}^{(j)},Y_{iK}^{(j)}) \right].
\end{multline}

\bibliographystyle{elsarticle-num}
\bibliography{references_JOTA}

\end{document}